\documentclass[12pt]{amsart}
\usepackage{amssymb}
\usepackage{amsthm}
\usepackage{amscd}

\usepackage{mathrsfs}
\usepackage{amsmath}
\usepackage{epic,eepic}
\usepackage{graphicx}
\DeclareMathAlphabet{\mathpzc}{OT1}{pzc}{m}{it}

\usepackage{latexsym}

\usepackage{pifont}

\theoremstyle{plain}
\newtheorem{lemma}{Lemma}[section]
\newtheorem{prop}[lemma]{Proposition}
\newtheorem{thm}[lemma]{Theorem}
\newtheorem{cor}[lemma]{Corollary}
\newtheorem{aplemma}{Lemma~A.\hspace{-1.5mm}}
\newtheorem{approp}{Proposition~A.\hspace{-1.5mm}}
\newtheorem{apthm}{Theorem~A.\hspace{-1.5mm}}
\newtheorem{apcor}{Corollary~A.\hspace{-1.5mm}}
\newtheorem{intthm}{Theorem}

\usepackage[all]{xy}

\newtheorem{conj}[lemma]{Conjecture}

\theoremstyle{definition}

\newtheorem{rema}[lemma]{Remark}

\newtheorem{remb}{Remark}

\newtheorem{defi}[lemma]{Definition}
\newtheorem{exa}[lemma]{Example}
\newtheorem{aprem}{Remark~A.\hspace{-1.5mm}}
\newtheorem{apdefi}{Definition~A.\hspace{-1.5mm}}
\newcommand{\bde}{\begin{defi}}
\newcommand{\ede}{\end{defi}\vspace{1mm}}
\newcommand{\ble}{\begin{lemma}}
\newcommand{\ele}{\end{lemma}}
\newcommand{\bpr}{\begin{prop}}
\newcommand{\epr}{\end{prop}}
\newcommand{\bt}{\begin{thm}}
\newcommand{\et}{\end{thm}}
\newcommand{\bco}{\begin{cor}}
\newcommand{\eco}{\end{cor}}
\newcommand{\bre}{\begin{rema}}
\newcommand{\ere}{\end{rema}}
\newcommand{\brea}{\begin{rema}}
\newcommand{\erea}{\end{rema}\vspace{1mm}}
\newcommand{\breb}{\begin{remb}}
\newcommand{\ereb}{\end{remb}\vspace{1mm}}
\newcommand{\bex}{\begin{exa}}
\newcommand{\eex}{\end{exa}}
\newcommand{\bpf}{\begin{proof}}
\newcommand{\epf}{\end{proof}\vspace{1mm}}

\newcommand{\bade}{\begin{apdefi}}
\newcommand{\eade}{\end{apdefi}}
\newcommand{\bale}{\begin{aplemma}}
\newcommand{\eale}{\end{aplemma}}
\newcommand{\bapr}{\begin{approp}}
\newcommand{\eapr}{\end{approp}}
\newcommand{\bat}{\begin{apthm}}
\newcommand{\eat}{\end{apthm}}
\newcommand{\baco}{\begin{apcor}}
\newcommand{\eaco}{\end{apcor}}
\newcommand{\bare}{\begin{aprem}}
\newcommand{\eare}{\end{aprem}}

\newcommand{\be}{\begin{enumerate}}
\newcommand{\ee}{\end{enumerate}}
\newcommand{\bcd}{\[\begin{CD}}
\newcommand{\ecd}{\end{CD}\]}
\newcommand{\bit}{\begin{itemize}}
\newcommand{\eit}{\end{itemize}}
\newcommand{\bq}{\begin{quote}}
\newcommand{\eq}{\end{quote}}
\newcommand{\ba}{\begin{array}}
\newcommand{\ea}{\end{array}}
\newcommand{\mcA}{\mathcal{A}}
\newcommand{\mcB}{\mathcal{B}}

\newcommand{\mcD}{\mathcal{D}}

\newcommand{\mcF}{\mathcal{F}}

\newcommand{\mcH}{\mathcal{H}}

\newcommand{\mcM}{\mathcal{M}}

\newcommand{\mcO}{\mathcal{O}}

\newcommand{\mcS}{\mathcal{S}}
\newcommand{\mcT}{\mathcal{T}}

\newcommand{\mbA}{\mathbb{A}}

\newcommand{\mbD}{\mathbb{D}}

\newcommand{\mbF}{\mathbb{F}}

\newcommand{\mbP}{\mathbb{P}}
\newcommand{\mbQ}{\mathbb{Q}}
\newcommand{\mbR}{\mathbb{R}}

\newcommand{\mbZ}{\mathbb{Z}}

\newcommand{\mfX}{\mathfrak{X}}
\newcommand{\mfY}{\mathfrak{Y}}
\newcommand{\mfZ}{\mathfrak{Z}}

\newcommand{\mpf}{\mathpzc{f}}

\newcommand{\msE}{\mathscr{E}}
\newcommand{\msF}{\mathscr{F}}
\newcommand{\msG}{\mathscr{G}}

\newcommand{\msK}{\mathscr{K}}
\newcommand{\msL}{\mathscr{L}}
\newcommand{\msM}{\mathscr{M}}
\newcommand{\msN}{\mathscr{N}}
\newcommand{\msO}{\mathscr{O}}

\newcommand{\msV}{\mathscr{V}}

\newcommand{\migi}{\rightarrow}
\newcommand{\longmigi}{\longrightarrow}

\newcommand{\hidari}{\leftarrow}

\newcommand{\isom}{\stackrel{\sim}{\migi}}

\newcommand{\migiincl}{\hookrightarrow}

\newcommand{\migisurj}{\twoheadrightarrow}

\newcommand{\mr}{\mathrm}
\newcommand{\hidden}[1]{\,}

\newcommand{\SSP}{\vspace{0mm}}
\newcommand{\LSP}{\vspace{0mm}}

\begin{document}

\title[Holonomic $\mathcal{D}$-modules of arithmetic type and middle convolution]{Holonomic $\mathcal{D}$-modules of arithmetic type \\ and middle convolution}
\author{Yasuhiro Wakabayashi}
\date{}
\markboth{Yasuhiro Wakabayashi}{}
\maketitle
\footnotetext{Y. Wakabayashi: 
Graduate School of Information Science and Technology, Osaka University, Suita, Osaka 565-0871, Japan;}
\footnotetext{e-mail: {\tt wakabayashi@ist.osaka-u.ac.jp};}
\begin{abstract}

The aim of the present paper is to study arithmetic properties of $\mathcal{D}$-modules on an algebraic variety over the field of algebraic numbers. We first provide a framework for extending a class of $G$-connections (resp., globally nilpotent connections; resp., almost everywhere nilpotent connections) to holonomic $\mathcal{D}$-modules. It is shown that the derived category of $\mathcal{D}$-modules in each of such extended classes carries a Grothendieck six-functor formalism. This fact leads us to obtain the stability of the middle convolution for $G$-connections with respect to the global inverse radii. As a consequence of our study of middle convolution, we prove equivalences between  various arithmetic properties on rigid Fuchsian systems. This result gives, for such systems of differential equations, an affirmative answer to a conjecture described in a paper written by  Y. Andr\'{e} and F. Baldassarri.

\end{abstract}
\tableofcontents 

\section{Introduction}
\LSP

\subsection{What is a $G$-connection?} \label{S01}

In the present paper, we  investigate
$G$-connections, i.e., certain connections on vector bundles satisfying a condition of moderate growth,
 that are fundamental in understanding the diophantine properties of 
   Siegel's $G$-functions.

Let $K$ be a number field, and denote by $K (t)$ the rational function field in one variable $t$  over $K$.
 Each $n \times n$ matrix  
$A \in M_n (K (t))$  (where $n \in \mbZ_{>0}$) with entries in $K(t)$ associates
a system of linear differential equations
\begin{align} \label{eQ21}
\mpf : \frac{d}{dt} \vec{y} = A \vec{y}, \hspace{5mm} \vec{y} = \begin{pmatrix}  y_1 \\ y_2 \\ \vdots \\ y_n\end{pmatrix}.
\end{align}
Solutions to this system may be identified with
horizontal elements (i.e., elements in the kernel) of
the connection 
 given by
\begin{align}
\nabla := \frac{d}{dt} - A : K(t)^n \migi K(t)^n. 
\end{align}
For each $s \in \mbZ_{\geq 0}$, we define $A_{[s]}$ as the matrix such that,
 if $\vec{y}$ is a solution of the system $\mpf$,
   then the equality $\left(\frac{d}{dt}\right)^s\vec{y} = A_{[s]} \vec{y}$ holds.
Hence,  $A_0$ coincides with the identity matrix $I_n$ and $A_{[s]}$ ($s =1,2, \cdots$) satisfy the recurrence relation: 
\begin{align}
A_{[s+1]} = \frac{d}{dt}A_{[s]}+ A \cdot A_{[s]}.
\end{align}

Given a finite place $v$ of $K$, we denote by $|-|_v$ the normalized non-archimedean absolute value corresponding to $v$.
The Gauss norm  $|- |_{\mr{Gauss}, v}$ on $K (t)$ determined by   $|-|_v$ (cf. (\ref{dd00})) extends to the norm on $M_n (K(t))$, which we also  denote by 
$\left|\!\left| - \right|\! \right|_{\mr{Gauss}, v}$.
Then, the {\bf global inverse radius} of $\nabla$
  is defined as 
\begin{align}
\rho (\nabla)  
:=  \sum_{v} \mr{log} \left(\mr{max} \left\{1, \limsup_{s  \to \infty} \left|\!\left|\frac{A_{[s]}}{s!}\right| \!\right|_{\mr{Gauss}, v}^{\frac{1}{s}} \right\} \right) \in \mbR_{\geq 0} \sqcup \{ \infty \}
\end{align}
(cf. (\ref{Eq33}) for a general definition), where the sum in the right-hand side runs over 
 the set of finite places $v$ of $K$.
We say that  $\nabla$ is a {\bf $G$-connection} (or, a {\bf $G$-operator})
if $\rho (\nabla) < \infty$.

Since the value $\rho (\nabla)$ is invariant under any base-change  over $K$, 
the notion of a $G$-connection can be formulated for connections on $\overline{\mbQ}(t)$-vector spaces (where $\overline{\mbQ}$ denotes the field of algebraic numbers) in a well-defined manner; moreover, we can generalize that notion 
to connections on a smooth algebraic variety over $\overline{\mbQ}$.
The class of $G$-connections has various basic examples.
In fact, 
it  contains 
 differential operators
 of minimal order annihilating   $G$-functions, e.g., 
  algebraic functions over $\mbQ (t)$ regular at the origin, the polylogarithm functions, and
 some hypergeometric series with rational parameters.
 Also, $G$-connections have specific properties, and they are 
  believed to come from geometry (Bombieri-Dwork's conjecture).
The study of $G$-connections with geometric treatments was 
  developed 
 from 1970s with the works of Galochkin, Chudnovsky, Andr\'{e},  Dwork, Baldassarri and other mathematician (cf. e.g., ~\cite{And1}, ~\cite{AB}, ~\cite{DGS}).

On the other hand, there are  several  other classes of connections characterized by 
 important arithmetic conditions, e.g.,
  {\it globally nilpotent connections} and {\it almost everywhere (a.e.)\,nilpotent connections} (cf. Definition \ref{Def1}).

\LSP
\subsection{First result: Generalization to holonomic $\mcD$-modules} \label{S011}

The primary purpose of the present paper is to 
 generalize $G$-connections, as well as globally nilpotent or a.e.\,nilpotent connections,  to holonomic $\mcD$-modules in order  to discuss various  functors (including the middle convolution functor) between the derived categories  of $\mcD$-modules from an arithmetic point of  view.

Given a smooth algebraic variety $X$ over   $\overline{\mbQ}$,
we denote by $D^b_h (\mcD_X)$ the derived category of bounded chain complexes of  $\mcD_X$-modules
having holonomic cohomology.
If $f : X \migi Y$ is a morphism of smooth varieties over $\overline{\mbQ}$, then we can push-forward $\mcD_X$-modules, as well as   pull-back $\mcD_Y$-modules,  along that morphism.
There are also some other natural functors between categories of $\mcD$-modules which together make up a version of the so-called ``six-functor formalism" of Grothendieck.

In the present paper, we introduce the full subcategory 
 \begin{align}
 D^b_{h} (\mcD_X)_G  \ \left(\text{resp.}, \ D^b_{h}(\mcD_X)_{\mr{nilp}}; \text{resp.,} \ D^b_{h} (\mcD_X)_{\mr{aen}}  \right)
 \end{align}
(cf. (\ref{Eq2201}))
 of $D^b_h (\mcD_X)$ consisting of bounded chain complexes
 such that $X$ can be stratified by locally closed subschemes on each of which the 
   cohomology sheaves are 
   $G$-connections (resp., globally nilpotent connections; resp., a.e.\,nilpotent connections).

Here, recall a result by Andr\'{e} and Baldassarri (cf. ~\cite[Main Theorem]{AB}), asserting  that
the cohomology sheaves of the push-forward of a $G$-connection define  {\it generically} $G$-connections.
Our formulation  in terms of $\mcD$-modules has the advantage that the result by    Andr\'{e}-Baldassari can be simply  formulated as the stability of the subcategories $ D^b_{h} (\mcD_{X})_G$ with respect to the push-forward functor. 
Together with other kinds of functors, we obtain the following assertion, which is the main result of the first part.

\SSP
\begin{intthm} [cf.  Theorems \ref{Th40} and \ref{Thr4}] \label{Theer}
Let $f : X \migi Y$ be a morphism of smooth algebraic varieties over $\overline{\mbQ}$.
Then,
there  is a six-functor formalism
\begin{align}
\int_f &: D_{h}^b (X)_G \migi D_{h}^b (Y)_G,  \\
  f^\dagger &:  D_{h}^b (Y)_G \migi  D_{h}^b (X)_G,  \notag\\
\int_{f !} &: D_{h}^b (X)_G \migi D_{h}^b (Y)_G,  \notag 
\end{align}
\begin{align}
  \ \ \ \ \ \ \ \ \ \ \     \ \  \ f^!    &:  D_{h}^b (Y)_G \migi  D_{h}^b (X)_G,  \notag \\
\mbD &:  D_{h}^b (X)_G^{\mr{op}}  \migi D_{h}^b (X)_G,\notag \\
\otimes^L_{\mcO_X} &: D_{h}^b (X)_G \times D_{h}^b (X)_G \migi D_{h}^b (X)_G \notag 
\end{align}
satisfying all the usual adjointness properties that one has in the theory of the derived category of $\mcD$-modules.
Also, the same assertion holds for $D^b_{h}(\mcD_X)_\mr{nilp}$ and $D^b_h (\mcD_X)_{\mr{aen}}$.
 \end{intthm}

\LSP
\subsection{Second result: Middle convolution on $G$-connections} \label{S015}

The second part of the present paper discusses
the middle convolution functors on the derived  categories under consideration.
The middle convolution, which is introduced by Katz (cf. ~\cite{Kat2}), is an operation for local systems on (an open subscheme of) the affine line and plays a fundamental role in the theory of rigid local systems.
By the Riemann-Hilbert correspondence, it can be formulated as  an operation on flat bundles, 
which
moreover carries 
a chain complex of $\mcD$-modules
 to another one.
The middle convolution depends on a parameter $\lambda$, and we denote by
\begin{align}
\mr{mc}_\lambda ( \msF^\bullet)
\end{align}
(cf. (\ref{Eq400k})) the result of that operation applied to a chain  complex $\msF^\bullet \in D^b_h (\mcD_{\mbA^1})$, where $\mbA^1$ denotes the affine line over $\overline{\mbQ}$.

As a corollary of Theorem \ref{Theer},  it is shown that
the assignment $\msF^\bullet \mapsto \mr{mc}_\lambda (\msF^\bullet)$ preserves the  arithmetic properties mentioned above.
In ~\cite{DR2},
Dettweiler and 
Reiter 
provided an explicit description of that operation on a  Fuchsian system  and  then
 studied
 how its $p$-curvature (for a prime $p$)
    changes under the convolution process.
 The results of that work enables us 
 to 
measure
 the complexity of the middle convolution for globally nilpotent or a.e.\,nilpotent connections   from an arithmetic point of view (cf.  ~\cite[Theorem 1]{DR2}).

On the other hand, to investigate  $G$-connections, we estimate the global inverse radii  $\rho (-)$ of  connections (cf. (\ref{Eq33}), (\ref{dd32})) by applying some basic generalities on  $\mcD$-modules together with the work of   Andr\'{e} and Baldassarri (cf. ~\cite[Theorem 3.1.2]{AB}).
Our results are summarized as follows.
 
\SSP
\begin{intthm} [cf. Theorems \ref{T491} and \ref{TT5}] \label{ThCC}
Let $\lambda$ be an element of $\mbQ \setminus \mbZ$.
\begin{itemize}
\item[(i)]
Let $\Box \in \left\{ G, \mr{nilp}, \mr{aen}\right\}$.
Then, 
the endofunctor $\mr{mc}_\lambda$ on $D_h^b (\mcD_{\mbA^1})$ given by assigning $\msF^\bullet \mapsto \mr{mc}_\lambda (\msF^\bullet)$ restricts to an endofunctor
\begin{align} \label{Eq991}
\mr{mc}_\lambda : D^b_h (\mcD_{\mbA^1})_\Box \migi   D^b_h (\mcD_{\mbA^1})_\Box
\end{align}
on $D_h^b (\mcD_{\mbA^1})_\Box$.
\item[(ii)]
Let $\msF$ be a flat bundle on a nonempty open subscheme $U$ of $\mbA^1$ of rank $n \in \mbZ_{>0}$.
Denote by $\iota$ the natural open immersion $U \migiincl \mbA^1$, and suppose that $\msF$ is of type $G$.
Then, the following inequalities of global inverse radii  hold: 
\begin{align}
\left| \rho (\mr{mc}_\lambda (\int_\iota \msF) ) -\rho (\mr{mc}_\lambda (\iota_{!*}\msF))   \right| \leq H (\lambda), \hspace{5mm} 
\rho (\mr{mc}_\lambda (\int_\iota \msF) ) \leq (n^2 +1) \left(\rho ( \msF) + H (\lambda)\right),
\end{align}
where $H (\lambda)$ denotes  the positive real number  defined in (\ref{df3}) and $\iota_{!*} \msF$ denotes the minimal extension of $\msF$ (cf. (\ref{Eq441})).
\end{itemize}
 \end{intthm}

\LSP
\subsection{Third result: Equivalence for rigid $G$-connections} \label{S0g11}

At the end of the present paper, we discuss an application of our study of  middle convolution.
Recall the following  conjecture, which  was  described in a paper by Andr\'{e} and  Baldassarri; 
 it asserts  an equivalence among
 the classes  of $G$-connections, globally nilpotent connection, and a.e.\,nilpotent connections.

\SSP
\begin{conj}[cf. ~\cite{AB}, \S\,1.4] \label{Conj1}
Let  $X$ be  a smooth algebraic variety over $\overline{\mbQ}$ and 
$\msF$ a flat  bundle on $X$.
Then, the following three conditions are equivalent to each other:
\begin{itemize}
\item[(a)]
$\msF$ is a $G$-connection;
\item[(b)]
$\msF$ is globally nilpotent;
\item[(c)]
$\msF$ is almost everywhere nilpotent.
\end{itemize}
 \end{conj}
\SSP

    To consider 
  this problem, 
 we focus on  {\it rigid} $G$-connections defined on an open subscheme $U$ of the projective line $\mbP^1$ over $\overline{\mbQ}$.
 (A connection is rigid if it is determined up to isomorphism by the conjugacy classes of its local monodromies; see \S\,\ref{SS028}.)
Katz 
 showed that any rigid irreducible local system on $U$
can be obtained from a rank one connection by applying iteratively a suitable sequence of middle convolutions and scalar multiplications 
  (cf. ~\cite{Kat2}, ~\cite{Ari}).
By applying this fact and 
the stabilities of the arithmetic classes of $\mcD$-modules proved in Theorem \ref{ThCC},
  we  give an   affirmative answer to  the above conjecture for rigid flat bundles, as described below. 
 In particular, 
 this provides a characterization of $G$-connections in terms of  $p$-curvature, i.e.,  
 without 
  using 
  any limits nor infinite sums.

\SSP
\begin{intthm} [cf.  Theorem \ref{TT90}] \label{ThM}
Conjecture \ref{Conj1} is true when $X$ is 
 a nonempty open subscheme $U$ of the projective line $\mbP^1$ over $\overline{\mbQ}$ and 
$\msF$ is a rigid flat bundle  on $U$.
 \end{intthm}

\vspace{10mm}
\section{Holonomic $\mcD$-modules of type $\mcA$}  \label{S1129}
\LSP

In this section, we discuss a formal  way  of generalizing   flat bundles 
satisfying suitable  conditions to chain complexes of holonomic $\mcD$-modules.
Given a subcategory  of the stack of flat bundles (i.e., an  ``$\mcA$" introduced in \S\,\ref{SS0233}),
we define  
the derived category of chain complexes of $\mcD$-modules such that the underlying variety  can be  stratified by locally closed subschemes on each of which the cohomology sheaves belong to that subcategory.
We also prove several  basic properties of chain complexes in  that  derived category from the homological point of view. 

By an {\it  algebraic variety} over a field $k_0$,
we  mean a geometrically connected  scheme  of finite type over $k_0$.
 We shall fix  an algebraically closed field $k$ of characteristic $0$.

\LSP
\subsection{Smooth stratifications} \label{SS02}

Let  $X$ be a nonempty connected reduced scheme of finite type  over $k$.

\SSP
\bde \label{Def2}
\begin{itemize}
\item[(i)]
A {\bf stratification} on $X$ is a collection $\mfX := \{ X_j \}_{j=0}^{m+1}$ (where $m \in \mbZ_{\geq 0}$) forming a decreasing sequence of reduced closed subschemes
\begin{align} \label{Eq25}
X = X_0 \supsetneq  X_1 \supsetneq X_2 \supsetneq \cdots \supsetneq X_m \supsetneq X_{m+1} = \emptyset 
\end{align}
 of $X$.
 If we are given a stratification on $X$  indicated, say, by  $\mfX := \{ X_j \}_{j=0}^{m+1}$, then it is 
 occasionally considered as a  collection   $\{ X_j \}_{j \in \mbZ}$ indexed by  the elements of $\mbZ$  by putting $X_j := X$ (resp., $X_j := \emptyset$) for
 $j < 0$ (resp., $j > m+1$).
  Also, for each 
   $j \in \mbZ$, we shall denote by 
 \begin{align} \label{Eq03030}
 \iota_{\mfX, j} : X_j \setminus X_{j+1} \migiincl X
 \end{align}
 the natural immersion.
\item[(ii)]
A stratification $\mfX := \{ X_j \}_{j=0}^{m+1}$ is called {\bf smooth} if, for each $j =0, \cdots, m$,
 the (nonempty) subscheme $X_j \setminus X_{j+1}$ of $X$ is smooth over $k$.
\end{itemize}
\ede

\SSP
\bde \label{Def21}
Let $\mfX : = \{ X_j \}_{j=1}^{m+1}$ and $\mfX' : = \{ X'_j \}_{j=1}^{m'+1}$ be  stratifications on $X$.
\begin{itemize}
\item[(i)]
We shall say that $\mfX$  {\bf is subordinate to $\mfX'$}
 if,
for each  $j' \in \{0, \cdots, m' \}$, there exists $j \in \{ 0, \cdots, m\}$ such that
the immersion 
$\iota_{\mfX', j'} : X'_{j'}\setminus X'_{j'+1} \migiincl X$
 factors through
$\iota_{\mfX, j} : X_{j}\setminus X_{j+1} \migiincl X$.
\item[(ii)]
We shall say that 
$\mfX$  {\bf is strictly subordinate to $\mfX'$}
 if, for each $j \in \{ 0, \cdots, m \}$, there exists  $j' \in \{0, \cdots, m' \}$ satisfying $X_j = X_{j'}$.
 (One may immediately verify that if $\mfX$ is strictly subordinate to $\mfX'$, then  it  is also subordinate to $\mfX'$.) 
\end{itemize}
\ede
\SSP

We shall prove the following two  basic properties on (smooth) stratifications.

\SSP
\bpr \label{Prop34}
Let $\mfX := \{ X_j \}_{j=0}^{m+1}$  be a stratification on $X$.
Then, there exists a smooth  stratification 
 on $X$ to which $\mfX$ is strictly subordinate.
In particular, there always exists a smooth stratification on $X$.
\epr
\begin{proof}
The latter assertion immediately  follows  from the former assertion by considering the case where the stratification $\mfX$ is taken as  $X = X_0 \supsetneq X_1 = \emptyset$.
Hence, it suffices to prove the former assertion.

Since 
$k$ is of  characteristic zero, 
each irreducible component of 
a reduced  subscheme of $X$
 has a dense open subscheme that is smooth over $k$.
By applying this fact successively, 
we see that, for each $j \in \{0, \cdots, m \}$,
there exists a smooth stratification 
 \begin{align}
 X_j \setminus X_{j+1} = Y_{j, 0} \supsetneq   \cdots  \supsetneq  Y_{j,  M_j} = \emptyset
 \end{align}
 (where $M_j \in \mbZ_{> 0}$) on $ X_j \setminus X_{j+1} \left(\neq \emptyset \right)$.
 We shall set $Y_{j, l}^+ := Y_{j, l} \cup X_{j+1}$ ($l = 0, \cdots, M_j$).
 Then, the decreasing sequence
 \begin{align}
 X = X_0  \left(= Y_{0, 0}^+ \right) \supsetneq Y_{0, 1}^+ 
 \supsetneq \cdots \supsetneq Y_{0, M_0}^+ 
 \left(= Y_{1,  0}^+ \right)  \supsetneq Y_{1,1}^+  \supsetneq   \cdots \supsetneq Y_{1, M_1}^+\supsetneq   \cdots \supsetneq  Y_{m,  M_m}^+ = \emptyset
 \end{align}
 defines a smooth stratification on $X$.
 Moreover, since $X_j = Y_{j,  0}^+$,
 $\mfX$ is strictly subordinate to  this smooth stratification. 
 This completes the proof of Proposition \ref{Prop34}.
 \end{proof}

\SSP
\bpr \label{Prop3}
Let $\mfX^1, \mfX^2, \cdots, \mfX^n$ (where $n \in \mbZ_{>0}$ and $\mfX^l := \{ X^l_j \}_{j=0}^{m_l +1}$)
be stratifications on $X$.
Then, there exists a smooth stratification $\mfX := \{ X_j \}_{j=0}^{m+1}$ on $X$
such that every $\mfX^l$ ($l=1, \cdots, n$) is subordinate to $\mfX$.
Moreover, 
  we can choose  $\mfX$ as such
  that $X_{m_1}^1 = X_j$ for some $j \in \{ 1, \cdots, m +1 \}$.
\epr
\begin{proof}
First, we shall consider the former assertion.
By induction on $n$,
 it suffices to consider the case of  $n=2$.
We shall set 
$l := \mr{min} \{ m_1, m_2 \}$ and $Y : =X^1_{m_1 - l}\cap X_{m_2 -l}^2$.
Let us  consider the following  decreasing sequence consisting of reduced subschemes of $Y$:
\begin{align}
Y = X^1_{m_1 - l}\cap X_{m_2 -l}^2
&\supseteq 
X^1_{m_l-l} \cap X^2_{m_2 - l +1}
\supseteq
X^1_{m_l-l+1} \cap X^2_{m_2 - l +1} \\
& \supseteq 
X^1_{m_l-l+1} \cap X^2_{m_2 - l +2}
\supseteq X^1_{m_l-l+2} \cap X^2_{m_2 - l +2} \notag  \\
& \supseteq
\cdots  \notag \\
&\supseteq
X^1_{m_1} \cap X^2_{m_2 +1}
\supseteq
X^1_{m_1 +1} \cap X^2_{m_2 +1} = \emptyset. \notag
\end{align}
By removing the duplicate constituents from this sequence, 
we obtain a  stratification $\mfY$ on $Y$.
Moreover, 
if $m_1 > m_2$ (resp., $m_2 > m_1$),
then
we add this stratification to the sequence $X^1_0 \supsetneq X^1_1 \supsetneq  \cdots \supsetneq X^1_{m_1 - l -1} $ (resp., $X^2_0 \supsetneq X^2_1 \supsetneq \cdots \supsetneq X^2_{m_2 - l -1}$);
we denote the resulting stratification on $X$ by $\mfX'$.
According to  Proposition \ref{Prop34},
there exists a smooth  stratification $\mfX$ on $X$ to which  $\mfX'$ is (strictly) subordinate.
By construction, both $\mfX^1$ and $\mfX^2$ are  subordinate to $\mfX$.
This completes the proof of the former assertion.

The latter assertion follows from the construction of $\mfX$ discussed above.
\end{proof}

\LSP
\subsection{Derived category of holonomic $\mcD$-modules} \label{SS02}

We now move on to  the discussion  on $\mcD$-modules.
 Much of the notation used in the present paper follows ~\cite{HTT}, and  we will refer to  
 that reference at many places in this text until the end of \S\,\ref{S1131}.
Although the discussions in {\it loc.\,cit.} only  deal with algebraic varieties over the field of complex numbers,
 the various results shown there based on a purely algebraic treatment of $\mcD$-modules remain true even when the base field is replaced with any arbitrary algebraically closed  field of characteristic $0$, e.g., the field of algebraic numbers $\overline{\mbQ}$.

Let $R$ be a commutative ring   and $X$ a smooth scheme over $\mr{Spec}(R)$ of constant relative dimension $d \in \mbZ_{> 0}$.
Denote by $\Omega_{X/R}$ the sheaf of $1$-forms on $X$ over $R$
 and by $\mcT_{X/R}$ its dual, i.e., the sheaf of vector fields on $X$ over  $R$.
Also, the canonical bundle $\omega_X$  of $X/R$   is   the line bundle defined as $\bigwedge^{d}\Omega_{X/R}$.

Recall that,  for an $\mcO_X$-module $\mcF$,
 an {\bf $R$-connection} on $\mcF$ is defined to be  an $R$-linear morphism
$\nabla : \mcF \migi \Omega_{X/R} \otimes \mcF$ satisfying $\nabla (av) = da \otimes v + a \cdot \nabla (v)$ for any local sections $a \in \mcO_X$ and $v \in \mcF$.
An $R$-connection is called {\bf flat} (or, {\bf integrable}) if it has vanishing curvature  (cf., e.g., ~\cite[\S\,1]{Kat}  for the definition of curvature).
By a {\bf (rank $n$) flat bundle} on $X/R$, we mean
 a pair $\msF := (\mcF, \nabla)$ consisting of a (rank $n$) vector bundle  $\mcF$ on $X$,  i.e., a locally free $\mcO_X$-module (of rank $n$),  and a flat $R$-connection $\nabla$ on $\mcF$.

Let us fix   a smooth algebraic variety $X$ over $k$.
Denote by $\mcD_X$ the sheaf of differential operators on $X/k$.
Each  $\mcD_X$-module is, by definition, given as
 an $\mcO_X$-module  together with 
 a left $\mcD_X$-action 
   extending its $\mcO_X$-module structure.
The class of flat bundles on $X/k$  coincides with   the class of $\mcD_X$-modules
  whose underlying sheaves  are  vector bundles (cf. ~\cite[Theorem 1.4.10]{HTT}).

Denote by $D^b_h (\mcD_X)$ the derived category of bounded chain complexes $\msF^\bullet$ 
of $\mcD_X$-modules having holonomic cohomology.
For each chain complex of $\mcD_X$-modules $\msF^\bullet$ and each $i \in \mbZ$,
the $i$-th cohomology sheaf of $\msF^\bullet$ will be denoted by $\mcH^i (\msF^\bullet)$.

Given another smooth algebraic variety $Y$ over $k$ and   a morphism $f  : X \migi Y$  over $k$, 
we 
obtain  the inverse image functor  
\begin{align} \label{Eq200}
Lf^* : D^b_h (\mcD_Y) \migi D^b_h (\mcD_X)
\end{align}
given by assigning $\msF^\bullet \mapsto D_{X \migi Y} \otimes^L_{f^{-1}D_Y}f^{-1}\msF^\bullet$, where $\mcD_{X \migi Y}$ denotes the $(\mcD_X, f^{-1}\mcD_Y)$-bimodule $\mcO_X \otimes_{f^{-1}\mcO_Y}f^{-1}\mcD_Y$.
It induces  the shifted inverse image functor
\begin{align} \label{Eq201}
f^\dagger := L f^* [\mr{dim}X - \mr{dim}Y] : D^b_h (\mcD_Y) \migi D^b_h (\mcD_X)
\end{align}
(cf. ~\cite[Theorem 3.2.3, (ii)]{HTT}).

\SSP
\begin{rema} \label{Rem39}
We here recall  two properties on the shifted  inverse image functor that will be used in the present paper.

\begin{itemize}
\item[(a)]
Let $f : X \migi Y$ be a morphism of smooth algebraic varieties over $k$.
Also, let 
 $\msF := (\mcF, \nabla)$ be  a flat bundle on $Y$; we regard it as 
an object of $D^b_h (\mcD_Y)$ concentrated at degree $0$.
Since the forgetful functor  $D_h^b (\mcD_{(-)}) \migi D^b (\mcO_{(-)})$ (where $D^b (\mcO_{(-)})$ denotes the derived category of bounded chain complexes of $\mcO_{(-)}$-modules) is exact, 
 we obtain   the following identification of $\mcO_X$-modules:
\begin{align}
\mcH^l (f^\dagger \msF^\bullet) =\begin{cases} 
f^*\mcF \left(=D_{X \migi Y} \otimes_{f^{-1}\mcD_Y}f^{-1}\msF\right) & \text{if $l = -\mr{dim}X + \mr{dim}Y$}; \\
0 & \text{if $l \neq  -\mr{dim}X + \mr{dim}Y$}.
\end{cases}
\end{align}
\item[(b)]
For each $\msF^\bullet \in D_h^{b} (\mcD_X)$, there exists a smooth stratification
$\mfX := \{ X_j\}_j$ on $X$ such that the cohomology sheaves $\mcH^l (\iota_{\mfX, j}^{\dagger}\msF^\bullet)$ are (possibly zero) flat bundles for all $l$ and $j$ (cf. ~\cite[Theorem 3.3.1]{HTT}).
\end{itemize}
\end{rema}
\SSP

By taking account of Remark \ref{Rem39}, (b), we make the following definition.

\SSP
\bde \label{Def3}
Let $\msF^\bullet$ be a chain complex in $D_h^b (\mcD_X)$ and  $\mfX := \{ \mfX_j \}_{j}$   a smooth stratification on $X$.
We shall say that 
$\mfX$  is a {\bf stratification for $\msF^\bullet$}
 if $\mcH^l (\iota_{\mfX, j}^\dagger \msF^\bullet)$ are flat bundles for all $l$ and $j$.
\ede
\SSP

Next,  
recall from ~\cite[Theorem 3.2.3, (i)]{HTT} that, for each $f: X \migi Y$ as above, 
we have  the direct image functor
\begin{align} \label{Ewgg3}
\int_f : D_h^b (\mcD_X) \migi D_h^b (\mcD_Y)
\end{align}
given by assigning
$\msF^\bullet \mapsto \int_f \msF^\bullet := R f_* (\mcD_{Y \hidari X} \otimes_{\mcD_X}^L \msF^\bullet)$, where 
$\mcD_{Y \hidari X}$ denotes the $(f^{-1}\mcD_Y, \mcD_X)$-module  $\omega_X \otimes_{\mcO_X} \otimes \mcD_{X \migi Y} \otimes_{f^{-1}\mcO_Y} f^{-1}\omega_Y^\vee$.

\LSP
\subsection{Holonomic $\mcD$-modules of type $\mcA$} \label{SS0233}

Denote by $\mcS m_{k}$ the category whose objects are smooth algebraic varieties over $k$ and whose morphisms are smooth $k$-morphisms between them.
Also, denote by $\mcF \mcB_k$ the category  over $\mcS m_{k}$ defined as follows:
\begin{itemize}
\item
The objects are the pairs $(Y, \msF)$ consisting of a smooth algebraic variety $Y$ over $k$ and a flat bundle $\msF$ on $Y$;
\item
The morphisms from $(Y, \msF)$ to $(Z, \msE)$ are the pairs $(f, \nu)$ consisting of a smooth morphism $f : Y \migi Z$ over $k$ and an isomorphism of flat bundles $\nu : \msF \isom f^*\msE$;
\item
The projection $\mcF \mcB_k \migi \mcS m_k$ is given by assigning $(Y, \msF) \mapsto Y$.
\end{itemize}
It is verified that $\mcF \mcB_k$ forms a category fibered in groupoids over $\mcS m_k$
by putting $(Y, f^*\msE)$ (for each morphism $f : Y \migi Z$ and $(Z, \msE) \in \mr{ob}(\mcF \mcB_k)$) as the {\it pull-back of $(Z, \msE)$  via $f$}.

 Now, 
let us fix  a   fibered subcategory  
\begin{align} \label{EEE33}
\mcA
\end{align}
   of $\mcF \mcB_k$
  having  the following three properties:
\begin{itemize}
\item[($\alpha$)]
Let $Y$ be an arbitrary  smooth algebraic variety over $k$.
Then, the category
$\mcA (Y)$ (i.e., the  fiber of the projection $\mcA \migi \mcS m_k$ over $Y$) is closed under taking flat subbundles, flat quotient bundles, and the duals of flat bundles.
Also,  $\mcA (Y)$  is closed under taking   extensions and  the tensor products of two flat bundles;
\item[($\beta$)]
If $f : Y \migi Z$ is a smooth morphism in $\mcS m_k$
 and $\msF$ is an object of $\mcA (Y)$ (considered as an object of $D^b_h  (\mcD_Y)$),
then there exists a dense open subscheme $U$ of $Z$ such that the cohomology sheaves $\mcH^l (\int_f \msF)|_U$ belong to $\mcA (U)$ for all $l$.
\item[($\gamma$)]
 If  $Y = \mr{Spec}(k)$, then  the equality  $\mcA (Y) = \mcF \mcB_k (Y)$ holds.
\end{itemize} 

\SSP
\begin{rema} \label{Rekki3}
The property that a flat bundle belongs to  $\mcA$ naturally extends to 
the situation where the underlying space is a disjoint union   $Y := \bigsqcup_{i=1}^r Y_i$ of smooth algebraic varieties $Y_i$.
To be precise,  we say, in the subsequent discussion,  that a flat bundle $\msF$ on such a scheme $Y$ belongs to $\mcA$ if the  restriction $\msF |_{Y_i}$ to each component $Y_i$ belongs to $\mcA (Y_i)$.
\end{rema}

\SSP
\bde \label{Def11}
\begin{itemize}
\item[(i)]
Let $\msF^\bullet$ be a chain  complex in $D_h^b (\mcD_X)$.
We shall say that $\msF^\bullet$  is  
 {\bf of type $\mcA$} 
if 
there exists a smooth stratification $\mfX := \{ X_j \}_j$  for $\msF^\bullet$
 such  that 
 the flat bundle $\mcH^l (\iota^\dagger_{\mfX, j} \msF^\bullet)$  belongs to $\mcA (X_l \setminus X_{l+1})$ for every  $l$ and $j$.
 In this situation,  $\mfX$ is called  {\bf an $\mcA$-stratification for  $\msF^\bullet$}.
(In particular, the notion of an $\mcF \mcB_k$-stratification, i.e., the case of $\mcA = \mcF \mcB_k$,  is nothing but  the notion of a stratification in the sense of Definition \ref{Def3}.)
\item[(ii)]
Let $\msF$ be a holonomic $\mcD_X$-module.
We shall say that $\msF$ is {\bf of type $\mcA$}
if it is of type $\mcA$ when considered as an object of $D^b_h (\mcD_X)$.
\end{itemize}
 \ede
\SSP

\begin{rema} \label{Rem42}
Let $Y$ be another smooth algebraic variety over $k$ and $f :X\migi Y$  a morphism  over $k$.
Also, let  $\msF := (\mcF, \nabla)$ be a flat bundle in  
 $\mcA (Y)$.
According to  Remark \ref{Rem39}, (a),
$f^\dagger \msF$ is quasi-isomorphic to
the inverse image  $f^* \msF$ shifted by $-\mr{dim}X+ \mr{dim}Y$.
Since the inverse image  functor $f^* (-)$ sends $\mcA (Y)$ to $\mcA (X)$,
$f^\dagger\msF$ belongs to $\mcA (X)$.
 By applying this fact to the case where $f = \iota_{\mfY, j}$ for a smooth stratification $\mfY := \{ Y_j \}_j$ on $Y$ (hence $X = Y_j \setminus Y_{j+1}$), we see that {\it any smooth stratification on $Y$  is an $\mcA$-stratification for  $\msF$}.
Moreover, this implies that 
{\it a flat bundle $\msF$ on $X$ belongs to $\mcA (X)$ if and only if
$\msF$  is of type $\mcA$ in the sense of  Definition \ref{Def11}}.
\end{rema}
\SSP

\bpr \label{Pro8}
Let $\msF^\bullet$ be  a chain complex 
in  $D^b_h(\mcD_X)$  and $\mfX := \{ X_j \}_{j=0}^{m+1}$, $\mfX' := \{X'_{j'}\}_{j'= 0}^{m'+1}$ two  smooth stratifications on $X$ such that $\mfX$ is subordinate to $\mfX'$.
Suppose that $\msF^\bullet$ is of type $\mcA$ and that $\mfX$ is
an $\mcA$-stratification for $\msF^\bullet$.
Then, 
 $\mfX'$  is an $\mcA$-stratification for $\msF^\bullet$.
 \epr
\begin{proof}
Let us take an arbitrary   $j' \in \{ 0, \cdots, m' \}$.
Since $\mfX$ is  subordinate to $\mfX'$, 
 there exists $j \in \{ 0, \cdots, m \}$ such that
the immersion $\iota_{\mfX', j'} : X_{j'} \setminus X_{j'+1} \migiincl X$ factors through $\iota_{\mfX, j} : X'_{j}\setminus X'_{j+1} \migiincl X$.
Denote by 
\begin{align}
\iota : X_{j'}\setminus X_{j'+1} \migiincl X_{j} \setminus X_{j+1}
\end{align}
 the resulting immersion (hence $\iota_{\mfX, j} \circ \iota = \iota_{\mfX', j'}$).
For each $l \in \mbZ$,
the $\mcO_{X_{j}\setminus X_{j+1}}$-module $\mcH^l (\iota^\dagger_{\mfX, j}\msF^\bullet)$ is locally free by assumption, so
we have
\begin{align} \label{Esw32}
\mcH^l (\iota^\dagger_{\mfX', j'}\msF^\bullet) \left(=\mcH^l ((\iota_{\mfX, j} \circ \iota)^\dagger\msF^\bullet)   \right)
 \cong \iota^\dagger \mcH^l (\iota^\dagger_{\mfX, j}\msF^\bullet).
\end{align}
Since
$\mcH^l (\iota^\dagger_{\mfX, j}\msF^\bullet)$ belongs to $\mcA (X_{j}\setminus X_{j+1})$ and 
 the inverse image functor $\iota^\dagger$ sends $\mcA (X_{j}\setminus X_{j+1})$ to $\mcA (X'_{j'} \setminus X'_{j'+1})$ (cf. Remark \ref{Rem42}), we have $\iota^\dagger \mcH^l (\iota^\dagger_{\mfX, j}\msF^\bullet) \in \mcA (X'_{j'} \setminus X'_{j'+1})$.
 Hence,  (\ref{Esw32})  implies $\mcH^l (\iota_{\mfX', {j'}}^\dagger \msF^\bullet) \in \mcA (X'_{j'} \setminus X'_{j'+1})$.
It follows that $\mfX$ is an $\mcA$-stratification for $\msF^\bullet$, and 
this completes  the proof of this proposition.
\end{proof}

\SSP
\bpr \label{P01} 
Let $\msE^\bullet$, $\msF^\bullet$, $\msG^\bullet$ are chain complexes in $D^b_h (\mcD_X)$, and suppose that there exists
a distinguished triangle
\begin{align} \label{Ewqq2}
\msE^\bullet \migi \msF^\bullet \migi \msG^\bullet \xrightarrow{+1}
\end{align}
in $D^b_h (\mcD_X)$.
Also, suppose that
two of $\msE^\bullet$, $\msF^\bullet$, and $\msG^\bullet$
are of type $\mcA$.
 Then,  the remaining one is of type $\mcA$.
 \epr
\begin{proof}
We only consider the case where $\msE^\bullet$ and $\msG^\bullet$ are supposed to be  of type $\mcA$ because  the proofs of the other cases are entirely similar.
By Propositions \ref{Prop3} and \ref{Pro8} applied to $\mcA = \mcF \mcB_k$,
 there exists a smooth stratification $\mfX := \{ X_j \}_{j=0}^{m+1}$ for the three chain complexes $\msE^\bullet$, $\msF^\bullet$, $\msG^\bullet$ that is  moreover an $\mcA$-stratification 
 for both $\msE^\bullet$ and $\msG^\bullet$.
 Let us choose  $j \in \{0, \cdots, m \}$.
 Then, we obtain from (\ref{Ewqq2}) an exact  sequence of  $\mcD_{X_j\setminus X_{j+1}}$-modules
\begin{align}
\cdots 
\migi \mcH^{l-1} (\iota_{\mfX, j}^\dagger \msG^\bullet)
\migi \mcH^l (\iota_{\mfX, j}^\dagger \msE^\bullet) 
\migi \mcH^l (\iota_{\mfX, j}^\dagger \msF^\bullet) 
\migi  \mcH^l (\iota_{\mfX, j}^\dagger \msG^\bullet) 
\migi
 \mcH^{l+1} (\iota_{\mfX, j}^\dagger \msE^\bullet) 
\migi
 \cdots. 
\end{align}
By the definition of $\mfX$,
every $\mcD_{X_j\setminus X_{j+1}}$-module in this sequence forms  a flat bundle.
Hence,
the assertion follows from the assumption that the category $\mcA (X_j\setminus X_{j+1})$
is closed under taking flat subbundles, flat quotient bundles, and extensions of two flat bundles.
\end{proof}
\SSP

\bco \label{Ew221}
Let $\alpha : \msE^\bullet \migi \msF^\bullet$ be a morphism in $D_h^b (\mcD_X)$, and suppose that both $\msE^\bullet$ and $\msF^\bullet$ are of type $\mcA$.
Then, there exists a chain complex $\msG^\bullet \in D_h^b (\mcD_X)$ of type $\mcA$ together with a distinguished triangle
\begin{align}
\msE^\bullet \xrightarrow{\alpha} \msF^\bullet \migi \msG^\bullet \xrightarrow{+1}.
\end{align}
\eco
\begin{proof}
The assertion follows from Proposition \ref{P01} and the fact that $D^b_h (\mcD_X)$ forms a triangulated category (cf. ~\cite[Corollary 3.1.4]{HTT}).
\end{proof}
\SSP

According to  Corollary \ref{Ew221} just  proved, one can define
the full triangulated subcategory
\begin{align} \label{Eq2}
D_{h}^b (\mcD_X)_\mcA
\end{align}
 of $D_h^b (\mcD_X)$ consisting of
chain complexes of type $\mcA$.

\vspace{10mm}
\section{Grothendieck's six functors for holonomic $\mcD$-modules of type $\mcA$}  \label{S1131}
\LSP

In this section,
we discuss   various functors on the derived categories of $\mcD$-modules and examine  their stability properties with respect to a fixed fibered subcategory $\mcA$ of $\mcF \mcB_k$.
  Theorem \ref{Theer} can be obtained by combining  the results of this section  applied to the case where $\mcA$ is taken to be one of the categories classifying  certain arithmetic flat bundles.

Let $k$, $X$, and $\mcA$ be as in the previous section.

\LSP
\subsection{Direct image functor} \label{SS023}

The first section deals with the direct image functor.
To begin with, let us prove the following proposition concerning closed or open immersions.

\SSP
\bpr \label{L011} 
The following assertions hold:
\begin{itemize}
\item[(i)]
Let $\iota : X \migiincl Y$ be a closed immersion between smooth algebraic  varieties over $k$.
Then, for each $\msF^\bullet \in D_{h}^b (\mcD_X)$,
we have
\begin{align}
 \msF^\bullet \in D_{h}^b (\mcD_X)_\mcA \ \Longleftrightarrow \ 
  \int_\iota \msF^\bullet \in D_{h}^b (\mcD_Y)_\mcA.
 \end{align}
 \item[(ii)]
 Let $\eta : X \migi Y$ be an open immersion between smooth algebraic varieties over $k$.
 Then, for each $\msF^\bullet \in D_{h}^b (\mcD_X)$,
we have
\begin{align}
 \msF^\bullet \in D_{h}^b (\mcD_X)_\mcA \ \Longleftrightarrow \ 
  \int_\eta \msF^\bullet \in D_{h}^b (\mcD_Y)_\mcA.
 \end{align}
 \end{itemize}
 In particular, the same equivalence assertion 
  holds for every (locally closed)  immersion  between smooth algebraic varieties $X \migiincl  Y$. 
\epr
\begin{proof}
We first consider the implication ``$\Leftarrow$" in  assertion (i).
Suppose that 
$\int_\iota \msF^\bullet \in D_{h}^b (\mcD_Y)_\mcA$.
This means that  there exists  a smooth $\mcA$-stratification
$\mfY^1$ for $\int_\iota \msF^\bullet$.
On the other hand, let us consider a smooth stratification $\mfY^2 := \{ Y^2_{j} \}_{j=0}^2$ on $Y$ determined by $Y^2_0 := Y$, $Y^2_1 := \mr{Im}(\iota)$, and  $Y^2_2 := \emptyset$.
By Proposition \ref{Prop3}, we can find a smooth stratification
$\mfY := \{ Y_{j'} \}_{j'=0}^{m+1}$ on $Y$ such that 
both $\mfY^1$ and $\mfY^2$ are  subordinate to $\mfY$ and $Y_{j_0} = \left( Y^2_1 = \right) \mr{Im}(\iota)$ for some $j_0 \in \{0, \cdots, m \}$.
It follows from  Proposition \ref{Pro8} that $\mfY$ is an $\mcA$-stratification for $\int_\iota \msF^\bullet$.
By putting
$X_j :=  \iota^{-1}(Y_{j + j_0})$ ($j= 0, \cdots, m-j_0 +1$), we obtain 
 a smooth stratification $\mfX := \{ X_j \}_{j=0}^{m-j_0 +1}$ on $X$.
Then,
for each $j = 0, \cdots, m-j_0$, we have
\begin{align} \label{Ew201}
\iota^\dagger_{\mfX, j}\msF^\bullet 
\cong
\iota_{\mfX, j}^\dagger \iota^\dagger \int_\iota \msF^\bullet
\cong
 (\iota \circ \iota_{\mfX,  j})^\dagger \int_\iota \msF^\bullet \cong \iota_{\mfY, j + j_0}^\dagger \int_\iota \msF^\bullet
\end{align}
under the identification $X_j \setminus X_{j+1} = Y_{j+j_0} \setminus Y_{j + j_0 +1}$ via $\iota$.
In particular,  the cohomology sheaves $\mcH^l (\iota^\dagger_{\mfX, j}\msF^\bullet) \left(\cong \mcH^l (\iota_{\mfY, j + j_0}^\dagger \int_\iota \msF^\bullet) \right)$ belong to $\mcA (X_j\setminus X_{j+1})$ for all $l$.
This implies $\msF^\bullet \in D_h^b (\mcD_X)_\mcA$, and  thus
 we have proved the implication ``$\Leftarrow$" in (i).

Next, we shall consider the implication ``$\Leftarrow$" in  assertion (ii).
Suppose that $\int_\eta \msF^\bullet \in D_{h}^b (\mcD_Y)_\mcA$, i.e., there exists an $\mcA$-stratification $\mfY := \{ Y_{j'} \}_{j'=0}^{m'+1}$ for $\int_\eta \msF^\bullet$.
The decreasing sequence
\begin{align}
\eta (X) = Y_0 \cap \eta (X) \supseteq Y_1 \cap \eta (X) \supseteq \cdots \supseteq Y_{m'+1} \cap \eta (X) = \emptyset
\end{align}
determines 
a smooth stratification $\mfX := \{ X_j \}_{j = 0}^{m+1}$ on $X$ after  removing the duplicate constituents and identifying $X$ with $\iota (X)$.
Let us take an arbitrary $j \in \{ 0, \cdots, m \}$.
By the definition of $\mfX$, there exists $j' \in \{0, \cdots, m' \}$ with $X_j = Y_{j'}\cap \eta (X)$ (via the identification $X = \iota (X)$).
Hence, for each $l\in \mbZ$,
we have
\begin{align}
\mcH^l \left(\iota_{\mfX, j}^\dagger \msF^\bullet \right) 
& \cong 
 \mcH^l \left(\iota_{\mfX, j}^\dagger  \eta^\dagger \int_\eta \msF^\bullet \right)  \\
&\cong 
\mcH^l \left((\eta \circ \iota_{\mfX, j})^\dagger \int_\eta \msF^\bullet 
\right) \notag \\
&\cong 
\mcH^l \left(\left(\iota_{\mfY, j'}^\dagger \int_\eta \msF^\bullet \right) |_{Y_{j'}\cap \eta (X)}\right)  \notag \\
& \cong 
\mcH^l \left( \iota_{\mfY, j'}^\dagger \int_\eta \msF^\bullet \right) |_{Y_{j'}\cap \eta (X)}. \notag
\end{align}
It follows that 
$\mcH^l \left(\iota_{\mfX, j}^\dagger \msF^\bullet \right)$ belongs to $\mcA (X_j \setminus X_{j+1})$, which implies $\msF^\bullet \in D_h^b (\mcD_X)_\mcA$.
Thus, we have finished the proof of the implication ``$\Leftarrow$" in   (ii).

Finally, we shall consider the inverse implication   ``$\Rightarrow$" in assertion (i) (resp., (ii)).
Suppose that $\msF^\bullet \in D_h^b (\mcD_X)_\mcA$, i.e., there exists an $\mcA$-stratification $\mfX := \{ X_j \}_{j=0}^{m+1}$ for $\msF^\bullet$.
Let $\mfY := \{ Y_j \}_{j' = 0}^{m+2}$ be the smooth stratification 
on $Y$ determined by the condition that  $Y_{0} := Y$ and  $Y_{j'} := \iota (X_{j'-1})$ if  $j' = 1, \cdots, m+2$ (resp., $Y_{m+2} := \emptyset$ and  $Y_{j'} := \eta (X_{j'}) \cup (Y\setminus X)$ if $j' = 0, \cdots, m+1$). 
Since $Y_0 \setminus Y_1 = Y\setminus X$ (resp., $Y_{m+1} \setminus Y_{m+2} = Y \setminus X$),  we have $\iota_{\mfY, 0}^\dagger \int_\iota \msF^\bullet =0$ (resp., $\iota_{\mfY, m+1}^\dagger \int_\eta \msF^\bullet =0$) (cf. ~\cite[Proposition 1.7.1, (ii)]{HTT}).
Moreover, for each $j' = 1, \cdots, m+1$ (resp., $j' = 0, \cdots, m$) and each $l \in \mbZ$,
it is immediately verified 
 that 
$\mcH^l (\iota_{\mfY, j'}^\dagger \int_\iota \msF^\bullet)$ 
(resp., $\mcH^l (\iota_{\mfY, j'}^\dagger \int_\eta \msF^\bullet)$ )
 is isomorphic to $\mcH^l (\iota_{\mfX, j' -1}^\dagger \msF^\bullet)$ (resp., $\mcH^l (\iota_{\mfX, j' }^\dagger \msF^\bullet)$) via the natural identification $X_{j' -1} \setminus X_{j'} = Y_{j'} \setminus Y_{j' +1}$ (resp., $X_{j'} \setminus X_{j'+1} = Y_{j'} \setminus Y_{j' +1}$).
This implies  that $\mcH^l (\iota_{\mfY, j'}^\dagger \int_\iota \msF^\bullet)$ (resp., $\mcH^l (\iota_{\mfY, j'}^\dagger \int_\eta \msF^\bullet)$) belongs to $\mcA (Y_{j'} \setminus Y_{j' +1})$, and hence,  we have
$\int_\iota \msF^\bullet \in D_h^b(\mcD_Y)_\mcA$ (resp., $\int_\eta \msF^\bullet \in D_h^b(\mcD_Y)_\mcA$).
This completes the proof of the implication  ``$\Rightarrow$".
\end{proof}
\SSP

\bco \label{Co223}
Let $Z$ be a closed subscheme of $X$ and $\msF^\bullet$ a chain complex  in  $D^b_h (\mcD_X)_\mcA$.
Then, the chain  complex $R \Gamma_Z (\msF^\bullet)$ (cf. ~\cite[\S\,1.7]{HTT}) belongs to $D_h^b (\mcD_X)_\mcA$.
\eco
\begin{proof}
Write $U := X \setminus Z$ and write $\eta : U \migi X$ for the natural open immersion.
Let us consider the distinguished triangle
\begin{align}
R \Gamma_Z (\msF^\bullet) \migi \msF^\bullet \migi \int_\eta \eta^\dagger \msF^\bullet \xrightarrow{+1}
\end{align}
(cf. ~\cite[Proposition 1.7.1, (i)]{HTT}).
By the assumption on $\msF^\bullet$ together with   Proposition  \ref{L011}, (ii),
$\int_\eta \eta^\dagger \msF^\bullet \left(=  \int_\eta (\msF^\bullet |_U) \right)$ belongs to $D_h^b (\mcD_X)_\mcA$.
Hence, the assertion follows from Proposition \ref{P01}.
\end{proof}
\SSP

\bco \label{Co22fe}
Let $Z$ be a reduced closed subscheme of $X$
defined as a disjoint union of finitely many closed points.
Also, let
  $\msF$ be a holonomic $\mcD_X$-module supported on $Z$.
  Then, 
  $\msF$ is of type $\mcA$.
 \eco
\begin{proof}
The assertion follows from Proposition \ref{L011}, (i), ~\cite[Corollary 1.6.2]{HTT} (i.e., Kashiwara's equivalence), and the assumption  ``$\mcA (\mr{Spec}(k)) = \mcF \mcB_k (\mr{Spec}(k))$" on $\mcA$ (i.e., the property ($\gamma$) described in \S\,\ref{SS0233}).
\end{proof}
\SSP

The following assertion will be used in the proof of Proposition \ref{T05} described below and also asserts   a special case of that proposition.

\SSP
\ble \label{Lemmi}
Let $f : X \migi Y$ be  a smooth morphism between smooth algebraic varieties over $k$.
Then, for each flat bundle $\msF$ in $\mcA (X)$,
the direct image  $\int_f \msF$ belongs to $D^b_h (\mcD_Y)_\mcA$.
\ele
\begin{proof}
Since $\msF$ belongs to $\mcA (X)$ and $f$ is smooth,
it follows from the assumption on $\mcA$ (i.e., the property ($\beta$) described in \S\ref{SS0233}) that
  there exists a dense open subscheme $U_Y$ of $Y$ satisfying 
 $\mcH^l (\int_f \msF)  |_{U_Y} \in \mcA (U_Y)$ for every $l$.
 Let $\eta_X$ (resp.,  $\eta_Y$)  denote the natural open immersion $f^{-1}(U_Y) \migiincl X$ (resp., $U_Y \migiincl Y$).
  Note that
 \begin{align} \label{Ew9090}
 \int_{\eta_Y}\eta_Y^\dagger\int_f \msF \cong \int_{\eta_Y} \int_{f_Y}\eta_X^\dagger\msF \cong  \int_f \int_{\eta_X} \eta_X^\dagger \msF,
 \end{align}
 where  $f_Y : f^{-1}(U_Y) \migi U_Y$ denotes  the smooth  morphism  obtained by restricting $f$.
 Since  $\eta_Y^\dagger\int_f \msF  \left(= \left(\int_f \msF\right) |_{U_Y} \right) \in D_h^b (\mcD_{U_Y})_\mcA$, it follows from Proposition  \ref{L011}, (ii), that
 the leftmost, hence also the rightmost, of    (\ref{Ew9090}) belongs to $D^b_h (\mcD_Y)_\mcA$.
If we set $Z := Y \setminus U_Y$, which is a reduced subscheme of $Y$, then 
the inverse image $f^{-1}(Z)$ specifies a reduced closed subscheme of $X$ with $f^{-1}(Z) = X \setminus f^{-1}(U_Y)$ and 
we obtain 
the distinguished triangle
 \begin{align}
 \int_f R\Gamma_{f^{-1}(Z)}\msF \migi \int_f \msF \migi \int_f \int_{\eta_X} \eta_X^\dagger \msF \xrightarrow{+1}
 \end{align}
 induced from  $R\Gamma_{f^{-1}(Z)}\msF \migi \msF \migi \int_{\eta_X} \eta_X^\dagger \msF \xrightarrow{+1}$.
 Hence, by Proposition \ref{P01} together with  the fact that $\int_f \int_{\eta_X} \eta_X^\dagger \msF \in D_h^b (\mcD_Y)_\mcA$, 
 the problem is  reduced to proving that $\int_f R\Gamma_{f^{-1}(Z)}\msF  \in D^b_h (\mcD_Y)_\mcA$. 
 
 Next, let us take  
 a dense open subscheme $U_Z$ of $Z$ that is smooth over $k$.
 We shall suppose  that
$W := Z \setminus U_Z$ is nonempty  (hence the inequality $\mr{dim}W < \mr{dim}Z$ holds  if $Z$ has positive dimension).
 Also, write $\eta_Z$ (resp., $\eta_W$) for the natural  immersion $f^{-1}(U_Z) \migiincl X$ (resp., the natural open immersion $X \setminus W \migiincl X$), which defines a morphism between smooth algebraic varieties.
Let  $f_Z$ denote the smooth morphism  $f^{-1}(U_Z) \migi U_Z$ obtained by restricting $f$.
 Since $\eta_Z^\dagger \msF \left(= \msF |_{f^{-1}(U_Z)} \right)\in \mcA (f^{-1}(U_Z))$, 
  there exists a dense open subscheme  of $U_Z$ on which  
 $\mcH^l (\int_{f_Z} \eta_Z^\dagger \msF)$ belongs to $\mcA$ for every $l$.
 Hence, by possibly replacing $U_Z$ with its open subscheme, we may assume that
 $\int_{f_Z} \eta_Z^\dagger \msF \in D_h^{b}(\mcD_{U_Z})_\mcA$. 
Let us  consider  the distinguished triangle
 \begin{align}
 R\Gamma_{f^{-1}(W)} \msF \migi R \Gamma_{f^{-1}(Z)} \msF  \migi \int_{\eta_Z}\eta_Z^\dagger \msF \xrightarrow{+1}
 \end{align}
induced from the fact that  $\int_{\eta_Z}\eta_Z^\dagger  \msF\cong
\int_{\eta_W} \eta_W^\dagger R \Gamma_{f^{-1}(Z)}\msF$ (cf. ~\cite[Corollary 1.6.2]{HTT}).
 It yields, via $\int_f$,  
a distinguished triangle
 \begin{align} \label{Ewidk}
 \int_f R \Gamma_{f^{-1}(W)} \msF \migi 
 \int_f R \Gamma_{f^{-1}(Z)} \msF \migi   \int_f  \int_{\eta_Z}\eta_Z^\dagger  \msF \xrightarrow{+1}.
 \end{align}
 If  $\iota$ denotes the natural immersion $U_Z \migiincl Y$, then we have 
 \begin{align} \label{EQw}
 \int_f\int_{\eta_Z}\eta_Z^\dagger  \msF \cong  \int_{\iota}\int_{f_Z} \eta_Z^\dagger \msF.
 \end{align}
 The right-hand side, hence also the left-hand side,  of (\ref{EQw}) belongs to $D_h^b (\mcD_Y)_\mcA$ because of Proposition \ref{L011}, (i) and (ii), together with the fact that $\int_{f_Z} \eta_Z^\dagger \msF \in D_h^b (\mcD_{U_Z})_\mcA$.
   Hence, by Proposition \ref{P01} and (\ref{Ewidk}), the problem is reduced to proving that $\int_f R \Gamma_{f^{-1}(W)}\msF \in D^b_h (\mcD_Y)_\mcA$. 
  
  By repeating the argument in the previous paragraph, 
  the problem is eventually reduced to proving  that $ \int_f R \Gamma_{f^{-1}(W)} \msF \in D^b_h (\mcD_Y)_\mcA$  in the case of  $\mr{dim}W = 0$.
  But, this is clear from 
  Corollary \ref{Co22fe}.
  We have finished  the proof of this lemma.
\end{proof}

\SSP
\bpr\label{T05} 
Let $f : X \migi Y$ be a morphism of smooth algebraic varieties over $k$.
Then, the  direct image functor $\int_f : D^b_h (\mcD_X) \migi D^b_h (\mcD_Y)$ 
(cf. (\ref{Ewgg3}))
restricts to a functor
\begin{align} \label{Ty1}
\int_f : D_{h}^b (\mcD_X)_\mcA \migi D_{h}^b (\mcD_Y)_\mcA.
\end{align}
 \epr
\begin{proof}
We prove this assertion by induction on the dimension $\mr{dim}X$.
The base step, i.e., the case of $\mr{dim}X = 0$, follows from Proposition \ref{L011}, (i).
In what follows, we shall consider the induction step.

Let $\msF^\bullet$ be a chain complex in $D_h^b (\mcD_{X})_\mcA$.
The problem is to show  that $\int_f \msF^\bullet \in D_h^b (\mcD_Y)_\mcA$.
By induction on  the cohomological length of $\msF^\bullet$,
we may assume that $\msF^\bullet = \msF$ for a holonomic $\mcD_X$-module $\msF$ of type $\mcA$.
Since $k$ has characteristic zero, there exists 
a dense open subscheme $U$ of the scheme-theoretic image $\mr{Im}f$ that is smooth over $k$.
Also, 
there exists  
 a dense open subscheme $V$ of $f^{-1}(U)$ such that 
 the morphism  $f_U : V \migi U$ obtained by restricting 
 $f$ is smooth.
 After possibly replacing $V$ with its open subscheme,
 we may assume that $\msF |_V$ is a flat bundle in $\mcA (V)$.
 Write $\eta$ (resp., $\iota$) for the natural open  immersion $V \migiincl X$ (resp., the natural immersion $U \migiincl Y$).  It follows from Proposition \ref{L011} and Lemma \ref{Lemmi}
 that $\int_\iota \int_{f_U} \eta^\dagger \msF \left(= \int_\iota \int_{f_U} \msF |_V \right) \in D_h^b (\mcD_Y)_\mcA$.
 Since  $\int_f \int_\eta \eta^\dagger \msF$ is quasi-isomorphic to
 $\int_\iota \int_{f_U} \eta^\dagger \msF$, it also belongs to  $D_h^b (\mcD_Y)_\mcA$.
 
 If $V = X$, then we have  already finished the proof.
Hence, it suffices to consider the case where   the reduced closed subscheme  $Z := X \setminus V$ of $X$ is nonempty (and $\mr{dim}Z < \mr{dim}X$).
By  the distinguished triangle 
\begin{align}
\int_f R \Gamma_Z \msF \migi \int_f \msF \migi \int_f \int_\eta \eta^\dagger \msF \xrightarrow{+1}
\end{align}
induced from $R \Gamma_Z \msF \migi \msF \migi \int_\eta \eta^\dagger \msF \xrightarrow{+1}$,
 the problem is reduced to proving  that 
$\int_f R \Gamma_Z \msF \in D_h^b (\mcD_Y)_\mcA$.

Now, let us take a dense open subscheme $U_Z$ of $Z$ that is smooth over $k$.
Write $W := Z \setminus U_Z$ (considered as a reduced closed subscheme of $X$) and write $\iota_{Z}$ for the natural immersion $U_Z \migiincl X$.
After possibly replacing $U_Z$ with its open subscheme,
we may assume that
the cohomology sheaves $\mcH^l (\iota^\dagger_Z \msF)$
are  flat bundles in $\mcA (U_Z)$ for all $l$.
(In fact, let $\mfX := \{ X_j \}_j$ be an $\mcA$-stratification for $\msF$.
There exists $j$ such that $X_j/X_{j+1}$ contains the generic point of $Z$, i.e., contains a dense open subscheme $U_Z$ of $Z$.
If $\overline{\iota}_Z$ denotes the natural immersion $U_Z \migiincl  X_j\setminus X_{j+1}$, then we have
 \begin{align}
 \mcH^l (\iota^\dagger_Z \msF) \cong \mcH^l (\overline{\iota}^\dagger_Z  \iota^\dagger_{\mfX, j}\msF) \cong \overline{\iota}^\dagger_Z \mcH^l (\iota^\dagger_{\mfX, j}\msF)\cong \overline{\iota}^*_Z \mcH^l (\iota^\dagger_{\mfX, j}\msF),
 \end{align}
where the second ``$\cong$" follows from Remark \ref{Rem39}, (a), together with the fact that the $\mcH^l (\iota^\dagger_{\mfX, j}\msF)$'s are flat bundles.
Hence, since   $\mcH^l (\iota^\dagger_{\mfX, j}\msF) \in \mcA (X_j \setminus X_{j+1})$,
we obtain  the claim, as desired.)
Let us consider   
the  distinguished triangle
\begin{align}
R \Gamma_W\msF \migi  R \Gamma_Z\msF \migi \int_{\iota_{Z}}\iota_{Z}^\dagger \msF \xrightarrow{+1} 
\end{align}
induced from the fact that   $\iota_{Z}^\dagger R \Gamma_Z \msF \cong  \iota_{Z}^\dagger\msF$
(cf. ~\cite[Corollary 1.6.2 or Proposition 1.7.1, (iii)]{HTT}).
It yields, via $\int_f$,
a distinguished triangle
\begin{align}
\int_f R \Gamma_W\msF \migi \int_f R \Gamma_Z \msF \migi 
 \int_{f\circ \iota_{Z}}\iota_{Z}^\dagger  \msF \xrightarrow{+1}. 
\end{align}
Since
$\iota_Z^\dagger \msF$ belongs to $D_h^b (\mcD_{U_Z})_\mcA$,  
the induction hypothesis implies 
$\int_{f \circ \iota_{Z}}\iota_{Z}^\dagger  \msF \in D_h^b (\mcD_Y)_\mcA$. 
 Hence, 
 the problem is reduced to proving that $\int_f R \Gamma_W \msF \in D^b_h (\mcD_Y)_\mcA$.

By repeating the argument in the previous paragraph, the problem is eventually reduced to proving that  $\int_f R \Gamma_W \msF \in D^b_h (\mcD_Y)_\mcA$ in the case of $\mr{dim}W=0$.
But, this is clear from
Corollary \ref{Co22fe}.
We have finished  the proof of this proposition.
\end{proof}

\LSP
\subsection{Inverse image functor} \label{SS023d}

Next, we prove that the inverse image functor preserves the subcategory  $\mcD_h^b (\mcD_X)_\mcA \subseteq D_h^b (\mcD_X)$.

\SSP
\bpr \label{Pr1} 
Let $f : X \migi Y$ be a morphism of smooth algebraic varieties over $k$.
Then, the  functor $f^\dagger$ (cf. (\ref{Eq201})) restricts to a functor
\begin{align}
f^\dagger : D^b_{h} (\mcD_Y) _\mcA\migi   D^b_{h} (\mcD_X)_\mcA.
\end{align}
Also, the same assertion holds for the functor $Lf^*$.
 \epr
\begin{proof}
By decomposing $f$ into the composite of the  closed immersion $(\mr{id}_X, f) :  X \migiincl X \times_k Y$ and the projection $X \times_k Y \migisurj Y$, we may assume
that $f$ is either a closed immersion or a smooth morphism.

We first consider the case where $f$ is smooth.
Let
$\msF^\bullet$ be a chain complex in   $D^b_{h} (\mcD_Y)_\mcA$, 
and choose an $\mcA$-stratification $\mfY := \{ Y_j \}_j$  for $\msF^\bullet$.
For each $j$, we denote by $f_j$ the morphism  $f^{-1}(Y_j) \migi Y_j$ obtained by restricting $f$.
Then,
the collection $f^{-1}\mfY := \{ f^{-1}(Y_j) \}_j$ forms a smooth stratification on $X$.
For each $l \in \mbZ$,
the smoothness of $f_j$ implies 
 that
\begin{align}
\mcH^l (\iota_{f^{-1}\mfY, j}^\dagger (f^\dagger \msF^\bullet)) \cong
\mcH^l (f_j^\dagger (\iota_{\mfY, j}^\dagger\msF^\bullet)) \cong
f_j^* \mcH^l (\iota_{\mfY, j}^\dagger \msF^\bullet).
\end{align}
In particular, 
the cohomology sheaves  $\mcH^l (\iota_{f^{-1}\mfY, j}^\dagger (f^\dagger \msF^\bullet))$ belong to $\mcA (f^{-1}(Y_j))$, so
$f^\dagger\msF^\bullet$ is verified to be an object of $D_h^b (\mcD_X)_\mcA$ (and $f^{-1}\mfY$ forms an $\mcA$-stratification for $f^\dagger\msF^\bullet$).

Next, let us consider the case of a closed immersion $\iota : X \migiincl Y$.
We shall write
 $\eta : U:= Y \setminus X \migiincl Y$ for  the  open immersion.
Let us consider 
the  distinguished triangle
\begin{align}
\int_\iota \iota^\dagger \msF^\bullet \migi \msF^\bullet \migi \int_\eta  \eta^\dagger \msF^\bullet \xrightarrow{+1}
\end{align}
(cf. ~\cite[Proposition 1.7.1, (i) and (iii)]{HTT}).
Since
 $\eta^\dagger \msF^\bullet \in D_h^b (\mcD_{Y \setminus X})_\mcA$ (by the previous argument for a smooth morphism $f$),
 Proposition  \ref{L011}, (ii),  implies 
$\int_\eta \eta^\dagger \msF^\bullet \in D^b_h (\mcD_Y)_\mcA$.
Hence, by 
 Proposition \ref{P01},
 $\int_\iota \iota^\dagger \msF^\bullet$ is verified to be a chain  complex in  $D^b_h (\mcD_Y)_\mcA$.
 By applying Proposition  \ref{L011}, (i), we have 
$\iota^\dagger \msF^\bullet \in D^b_h (\mcD_X)_\mcA$, thus completing 
 the proof of this proposition.
\end{proof}

\LSP
\subsection{Tensor product functor} \label{SS018}

Let $Y$ be another smooth algebraic variety over $k$.
Let 
$\pi_1$ (resp., $\pi_2$) denote the projection from the product of two $k$-schemes onto  
the first (resp., second) factor, e.g., the projection $X \times_k Y \migi X$ (resp., $X \times_k Y \migi Y$).

Recall that, for a $\mcD_X$-module $\msF$ and a $\mcD_Y$-module $\msE$, 
the {\it exterior  tensor product} of $\msF$ and $\msE$ is defined as the $\mcD_{X \times_k Y}$-module
\begin{align}
\msF \boxtimes \msE := \mcD_{X \times_k Y} \otimes_{\pi_1^{-1}\mcD_X \otimes_k \pi_2^{-1}\mcD_Y} (\pi_1^{-1}\msF \otimes_k \pi_2^{-1}\msE).
\end{align}
The assignment $(\msF, \msE) \mapsto \msF \boxtimes \msE$ extends to 
 a functor
\begin{align}
(-) \boxtimes (-) : D_{h}^b (\mcD_X) \times D_{h}^b (\mcD_Y) \migi D_{h}^b (\mcD_{X \times_k Y})
\end{align}
(cf. ~\cite[Proposition 3.2.2]{HTT}).

\SSP
\bpr \label{P04} 
The exterior tensor product $\boxtimes$
restricts to a functor
\begin{align}
 (-) \boxtimes (-) : D_{h}^b (\mcD_X)_\mcA \times D_{h}^b (\mcD_Y)_\mcA \migi D_{h}^b (\mcD_{X \times_k Y})_\mcA.
\end{align}
 \epr
\begin{proof}
Let $\msF^\bullet$ and $\msE^\bullet$ be chain complexes  in $D_h^b (\mcD_X)_\mcA$ and $D_h^b (\mcD_Y)_\mcA$, respectively.
Choose $\mcA$-stratifications $\mfX := \{ X_i \}_{i=1}^{m+1}$, $\mfY := \{ Y_j \}_{j=1}^{n+1}$ 
for
$\msF^\bullet$ and $\msE^\bullet$, respectively.
For each $j \in \{ 0, \cdots, m+n +1\}$,
we shall set  $Z_j$ to be  the reduced closed  subscheme of $X \times_k Y$ defined as the union
$\bigcup_{l=0}^j X_l \times_k Y_{j-l}$.
The subscheme $Z_j \setminus Z_{j+1}$ of $X \times_k Y$ decomposes as 
\begin{align}
Z_j \setminus Z_{j+1} = \bigsqcup_{l=0}^j (X_l \setminus X_{l+1}) \times_k (Y_{j-l} \setminus Y_{j-l+1}).
\end{align}
In particular,   
it is a disjoint union of smooth subschemes of $X \times_k Y$, and 
the resulting collection $\mfZ := \{ Z_j \}_{j=0}^{m+n+1}$ forms a smooth stratification on $X \times_k Y$. 

Now, let us take a connected component $U$ of $Z_j \setminus Z_{j+1}$, which coincides with 
$(X_s \setminus X_{s+1}) \times_k (Y_{j-s} \setminus Y_{j-s+1})$ for some $s$.
For each $a\in \mbZ$ (resp., $b \in \mbZ$),
the local freeness of
 $\mcH^a ( \iota^{\dagger}_{\mfX, s}\msF^\bullet)$ (resp., $\mcH^b (\iota^{\dagger}_{\mfY, j-s}\msE^\bullet)$)  implies that
 $\mcH^a (\pi_1^\dagger \iota^{\dagger}_{\mfX, s}\msF^\bullet)$
 (resp., $\mcH^b (\pi_2^\dagger \iota^{\dagger}_{\mfY, j-s}\msE^\bullet)$)
 is  isomorphic to
 $ \pi_1^*\mcH^a ( \iota^{\dagger}_{\mfX, s}\msF^\bullet) $ (resp., $\pi_2^*\mcH^b (\iota^{\dagger}_{\mfY, j-s}\msE^\bullet)$) and hence locally free.
Since the forgetful functor $D_h^b (\mcD_{(-)}) \migi D^b (\mcO_{(-)})$ is compatible with the exterior product functor $\boxtimes$ (cf. the discussion preceding ~\cite[Proposition 1.5.18]{HTT}),
 the natural morphism of $\mcD_{U}$-modules
 \begin{align}
 \mcH^l(\iota^{\dagger}_{\mfX, s}\msF^\bullet \boxtimes \iota^{\dagger}_{\mfY, j-s} \msE^\bullet) 
 \isom
  \bigoplus_{a+ b = l} \mcH^a (\pi_1^\dagger \iota^{\dagger}_{\mfX, s}\msF^\bullet) \otimes \mcH^b (\pi_2^\dagger \iota^{\dagger}_{\mfY, j-s}\msE^\bullet)
 \end{align} 
is an isomorphism because of the K\"{u}nneth formula for chain complexes. 
Thus, 
we obtain  the following  sequence of isomorphisms defined  for each $l$:
\begin{align}
\mcH^l (\iota_{\mfZ, j}^{\dagger} (\msF^\bullet\boxtimes \msE^\bullet)) |_U
&\isom  \mcH^l (\iota_{\mfZ, j}^{\dagger} (\msF^\bullet\boxtimes \msE^\bullet)|_U)
 \\
& \isom
 \mcH^l(\iota^{\dagger}_{\mfX, s}\msF^\bullet \boxtimes \iota^{\dagger}_{\mfY, j-s} \msE^\bullet) \notag \\
& \isom  
 \bigoplus_{a+ b = l} \mcH^a (\pi_1^\dagger \iota^{\dagger}_{\mfX, s}\msF^\bullet) \otimes \mcH^b (\pi_2^\dagger \iota^{\dagger}_{\mfY, j-s}\msE^\bullet) \notag\\
 & \isom
 \bigoplus_{a+ b = l} \pi_1^*\mcH^a ( \iota^{\dagger}_{\mfX, s}\msF^\bullet) \otimes \pi_2^*\mcH^b (\iota^{\dagger}_{\mfY, j-s}\msE^\bullet), \notag
\end{align}
where the second arrow follows from ~\cite[Proposition 1.5.18, (i)]{HTT}.
 Since both $\pi_1^*\mcH^a ( \iota^{\dagger}_{\mfX, s}\msF^\bullet)$ and $\pi_2^*\mcH^b (\iota^{\dagger}_{\mfY, j-s}\msE^\bullet)$
 belong to  $\mcA (U)$, we see that
 $\mcH^l (\iota_{\mfZ, j}^{\dagger} (\msF^\bullet\boxtimes \msE^\bullet)) |_U \in \mcA (U)$, which implies $\mcH^l (\iota_{\mfZ, j}^{\dagger} (\msF^\bullet\boxtimes \msE^\bullet)) \in \mcA (Z_j \setminus Z_{j+1})$.
 That is to say, 
 $\msF^\bullet\boxtimes \msE^\bullet$  lies in $D_h^b (\mcD_{X \times_k Y})_\mcA$ and $\mfZ$ forms an $\mcA$-stratification for $\msF^\bullet\boxtimes \msE^\bullet$.
This completes the proof of this proposition.
\end{proof}
\SSP

Also, the internal  tensor product $(\msF, \msE) \mapsto \msF \otimes_{\mcO_X} \msE$ induces 
the derived functor
\begin{align}
(-) \otimes_{\mcO_X}^L (-) : D^b_h (\mcD_X) \times D^b_h (\mcD_X) \migi D_h^b (\mcD_X).
\end{align}

\SSP
\bpr\label{C04} 
The functor  $\otimes_{\mcO_X}^L$ just mentioned  restricts to a functor
 \begin{align}
 (-) \otimes_{\mcO_X}^L (-) : D_{h}^b (\mcD_X)_\mcA \times D_{h}^{b} (\mcD_X)_\mcA \migi D_{h}^b (\mcD_X)_\mcA.
 \end{align}
 \epr
\begin{proof}
The assertion  follows from Propositions \ref{Pr1} and  \ref{P04}  because
$(-) \otimes^L_{\mcO_X} (-) = L \Delta_X^* ((-)\boxtimes (-))$, where $\Delta_X$ denotes the diagonal embedding $X \migiincl X \times X$.   
\end{proof}

\LSP
\subsection{Duality functor} \label{SS017}

Recall from ~\cite[Proposition 3.2.1]{HTT} 
the duality functor
\begin{align}
\mbD : D^b_{h} (\mcD_X) \isom D^b_{h} (\mcD_X)^{\mr{op}},
\end{align}
which is an equivalence of categories   given by assigning 
$\msF^\bullet \mapsto \mbD \msF^\bullet := R \mcH om_{\mcD_X} (\msF^\bullet, \mcD_X \otimes_{\mcO_X}\omega_X^{\vee} [\mr{dim}X])$.
If $\msF^\bullet = \msF$ for a flat bundle $\msF$ on $X$, then
$\mbD \msF$ may be identified  with the dual of $\msF$ in the usual sense.

\SSP
\bpr \label{P03} 
The  functor $\mbD$
restricts to an equivalence of categories
\begin{align}
\mbD : D^b_{h} (\mcD_X)_\mcA \isom D^b_{h} (\mcD_X)^{\mr{op}}_\mcA.
\end{align}
 \epr
\begin{proof}
Let $\msF^\bullet$ be a chain complex  in  $D^b_h (\mcD_X)_\mcA$, and 
choose an $\mcA$-stratification $\mfX := \{ X_j \}_{j=0}^{m+1}$ for  $\msF^\bullet$.
By descending   induction on $j$, we shall prove the claim  that
$\mbD R \Gamma_{X_j} \msF^\bullet$ lies in $D^b_h (\mcD_X)_\mcA$.

First, let us consider  the case of $j=m$, as the base step.
Since
$\mcH^{-l} (\iota^\dagger_{\mfX, m}\msF^\bullet)$ (for each $l$)  is a flat bundle in $\mcA (X_m)$,
its dual $\mbD \mcH^{-l} (\iota^\dagger_{\mfX, m}\msF^\bullet) \left(\cong \mcH^l (\mbD \iota^\dagger_{\mfX, m}\msF^\bullet) \right)$ belongs to $\mcA (X_m)$.
This implies  
$\mbD \iota^\dagger_{\mfX, m}\msF^\bullet\in D_h^b (\mcD_{X_m})_\mcA$.
Hence, by Proposition  \ref{L011},
$\int_{\iota_{\mfX, m}} \mbD\iota_{\mfX, m}^\dagger \msF^\bullet$ lies in 
$D^b_h (\mcD_X)_\mcA$.
On the other hand,  observe  that
\begin{align}
 \int_{\iota_{\mfX, m}} \mbD\iota_{\mfX, m}^\dagger \msF^\bullet
\cong
\mbD\int_{\iota_{\mfX, m}} \iota_{\mfX, m}^\dagger \msF^\bullet
\cong
\mbD R \Gamma_{X_m} \msF^\bullet, 
\end{align}
where the first ``$\cong$"
follows from ~\cite[Theorem 2.7.2]{HTT} together with the fact that $\iota_{\mfX, m}$ is a closed immersion, 
and the second ``$\cong$" follows from
the smoothness of $X_m$ together with  ~\cite[Proposition 1.7.1, (iii)]{HTT}.
It follows that 
$\mbD R \Gamma_{X_m} \msF^\bullet$ belongs to $D_h^b (\mcD_X)_\mcA$, which proves the base step.

Next, we shall consider the induction step.
To do this,  we suppose that we have proved the claim  with $j$ replaced by  $j+1$.
Let us consider the distinguished triangle
\begin{align}
R\Gamma_{X_{j+1}}\msF^\bullet \migi R \Gamma_{X_j}\msF^\bullet \migi  \int_{\iota_{\mfX, j}} \iota_{\mfX, j}^\dagger R \Gamma_{X_j}\msF^\bullet \xrightarrow{+1}.
\end{align}
By applying $\mbD (-)$ to it, we obtain
a distinguished triangle
\begin{align} \label{Eq22ko}
\mbD\int_{\iota_{\mfX, j}} \iota_{\mfX, j}^\dagger R \Gamma_{X_j}\msF^\bullet
 \migi 
 \mbD R \Gamma_{X_j}\msF^\bullet
  \migi
   \mbD R\Gamma_{X_{j+1}}\msF^\bullet    \xrightarrow{+1}.
\end{align}
Observe that
\begin{align} \label{Erood}
\mbD\int_{\iota_{\mfX, j}} \iota_{\mfX, j}^\dagger R \Gamma_{X_j}\msF^\bullet
\cong
\mbD \int_{\iota_\mfX, j} \iota^\dagger_{\mfX, j}\msF^\bullet
\cong
\int_{\iota_\mfX, j} \mbD\iota^\dagger_{\mfX, j}\msF^\bullet,
\end{align}
where the first ``$\cong$" follows from ~\cite[Proposition 1.7.1, (iii)]{HTT}.
Similarly to the argument in the base step,
we see that $\int_{\iota_\mfX, j} \mbD\iota^\dagger_{\mfX, j}\msF^\bullet$  belongs to
$D^b_h (\mcD_X)_\mcA$, so  (\ref{Erood}) implies 
$\mbD\int_{\iota_{\mfX, j}} \iota_{\mfX, j}^\dagger R \Gamma_{X_j}\msF^\bullet \in D^b_h (\mcD_X)_\mcA$.
Hence, by the induction hypothesis together with (\ref{Eq22ko}),
$\mbD R \Gamma_{X_j}\msF^\bullet$ is verified to be a chain complex in  $D^b_h (\mcD_X)_\mcA$.
The  induction step, hence also the claim, has been proved.

Consequently, we conclude that $\msF^\bullet \left(=R\Gamma_{X_0}\msF^\bullet \right)$ belongs to $D^b_h (\mcD_X)_\mcA$, and this completes the proof of this proposition.
\end{proof}

\LSP
\subsection{Proper direct/inverse image functors} \label{SS01g7}

Let $f  : X \migi Y$ be a  morphism of smooth algebraic varieties over $k$.
By using the duality functor $\mbD$, we obtain  functors
\begin{align}
\int_{f !} &:= \mbD \circ \int_f  \circ \mbD  : D_{h}^b (\mcD_X)_\mcA \migi D_{h}^b (\mcD_Y)_\mcA,
 \\
f^\star &:= \mbD  \circ f^\dagger \circ \mbD : D_{h}^b (\mcD_Y)_\mcA \migi D_{h}^b (\mcD_X)_\mcA, \notag
\end{align}
i.e., the proper direct and inverse image functors, respectively.
These functors together with $\int_f$ and $f^\dagger$ 
satisfy all the usual adjointness properties that one has in the theory of the derived category of $\mcD$-modules.
 In particular,
for each $\msF^\bullet \in D_{h}^b (\mcD_X)_\mcA$ and $\msE^\bullet \in D_{h}^b (\mcD_Y)_\mcA$,  there exist  natural isomorphisms
\begin{align}
\mr{Hom}_{D_h^b (\mcD_Y)_\mcA} 
 \left(\int_{f !} \msF^\bullet, \msE^\bullet \right) \cong \mr{Hom}_{D^b_h (\mcD_X)_\mcA} \left(\msF^\bullet, f^\dagger \msE^\bullet \right),
 \\
 \mr{Hom}_{D_h^b(\mcD_X)_\mcA} \left(f^\star \msE^\bullet, \msF^\bullet \right)
 \cong
 \mr{Hom}_{D^b_h (\mcD_Y)_\mcA} \left(\msE^\bullet, \int_{f}\msF^\bullet \right)\notag
\end{align}
(cf. ~\cite[Corollary 3.2.15]{HTT}).

As summarized below, we have now obtained  various  functors on the derived categories of $\mcD$-modules of type $\mcA$, namely
$\int_f$, $f^\dagger$, $\int_{f!}$, $f^!$, $\mbD$, $\otimes_{\mcO_X}^L$  (or $\boxtimes$), which forms 
 an example of the six-functor  formalism of Grothendieck (cf. ~\cite{Meb}).

\SSP
\bt \label{Th40}
Let  $\mcA$ be a fibered subcategory of $\mcF \mcB_k$ satisfying the three conditions ($\alpha$)-($\gamma$) described at the beginning of  \S\,\ref{SS0233}.
Then,
the full triangulated subcategory $D^b_h (\mcD_X)_\mcA$ of $D^b_h (\mcD_X)$ is stable under the functors $\int_f$, $f^\dagger$, $\int_{f!}$, $f^!$ (for each  morphism of smooth algebraic varieties $f : X \migi Y$), $\otimes^L_{\mcO_X}$, and $\mbD$.
\et

\LSP
\subsection{Minimal extensions} \label{SS037}

In this final subsection of \S\,\ref{S1131}, our discussion is restricted to 
 the case where {\it $X$ is a smooth curve, i.e., a smooth algebraic variety over $k$ of dimension $1$}.

Let $Y$ be a (locally closed) subvariety of $X$, and
denote by $\iota$ the natural immersion $Y \migiincl X$.
Also, let $\msF$  be  a holonomic $\mcD_Y$-module.
Since  $\mcH^l \left(\int_\iota \msF\right) = \mcH^l \left(\int_{\iota !} \msF\right) = 0$ for any $l \neq 0$,
both  $\int_\iota \msF$ and $\int_{\iota !} \msF$ can be regarded as $\mcD_X$-modules.
Recall that there exists a natural morphism
$\int_{\iota !} \msF \migi \int_\iota \msF$
(cf. ~\cite[Theorem 3.2.16]{HTT})
and  that
the {\it minimal extension} of $\msF$ is defined  as 
the image
\begin{align} \label{Eq441}
\iota_{!*}\msF := \mr{Im} \left( \int_{\iota !} \msF \migi \int_\iota \msF\right)
\end{align}
 of this morphism  (cf. ~\cite[Definition 3.4.1]{HTT}).

\SSP
\bpr \label{Prop909}
\begin{itemize}
\item[(i)]
Let  $Y$ and $\iota$ be as above and 
 and $\msF$ an irreducible  holonomic $\mcD_Y$-module of type $\mcA$.
Then, the minimal extension $\iota_{!*}\msF$ (which is irreducible, by ~\cite[Theorem 3.4.2]{HTT})  is of type $\mcA$.
\item[(ii)]
Any irreducible holonomic $\mcD_X$-module of type $\mcA$ is isomorphic 
to the minimal extension $\iota_{!*}\msF$ for some pair $(\iota, \msF)$, where
$\iota$
denotes the natural immersion $Y \migiincl X$ determined by 
 a locally closed subscheme $Y$ of $X$,  and  $\msF$ denotes an irreducible flat bundle on $Y$ of type $\mcA$.
\item[(iii)]
Let $\msF$ be a holonomic $\mcD_X$-module, and consider a finite sequence
\begin{align}
\msF = \msF_0 \supseteq \msF_1 \supseteq \cdots \supseteq \msF_m \supseteq \msF_{m+1} = 0
\end{align} 
of holonomic $\mcD_X$-submodules such that $\msF_j/\msF_{j+1}$ 
is irreducible for each $j$.
(Such a sequence always exists because of ~\cite[Proposition 3.1.2, (ii)]{HTT}.)
If $\msF$ is of type $\mcA$, then each $\msF_j/\msF_{j+1}$ is of type $\mcA$.
\end{itemize}
\epr
\begin{proof}
First, we shall consider assertion (i).
Note that $Y$ is either an open subscheme of $X$ or a disjoint union of finitely many closed points.
The latter  case follows immediately from 
Corollary \ref{Co22fe}.
In what follows, we shall suppose that $Y$ is an open subscheme of $X$.
Denote by $\msE$ the cokernel of the natural inclusion $\iota_{!*}\msF \migiincl \int_\iota \msF$.
In particular, we obtain a short exact sequence of holonomic $\mcD_X$-modules
\begin{align} \label{Eqrr3}
0 \longmigi \iota_{!*}\msF \longmigi \int_\iota \msF \longmigi \msE \longmigi 0.
\end{align}
Since $\iota^\dagger \iota_{!*}\msF = \iota^\dagger\int_\iota \msF = \msF$ and the functor $\iota^\dagger$ is exact,
we have $\iota^\dagger \msE = 0$.
It follows that
$\msE$ is supported in $X \setminus Y$ (which is a disjoint union of closed points).
According to  
Corollary \ref{Co22fe},
$\msE$ is of type  $\mcA$.
Hence, the exactness of (\ref{Eqrr3}) and Proposition \ref{P01} together imply 
 that $\iota_{!*}\msF$ is of type $\mcA$.
 This completes the proof of assertion (i).
 
Next, we shall consider assertion (ii).
Let $\msG$ be an irreducible  holonomic $\mcD_X$-module of type $\mcA$.
By ~\cite[Theorem 3.4.2, (ii)]{HTT}, it is isomorphic to $\iota_{!*}\msF$, where 
$\iota$ denotes the natural immersion $Y \migiincl X$ determined by
 a locally closed subscheme  $Y$ of $X$,  and $\msF$ denotes an irreducible  flat bundle on $Y$.
Moreover, since   $\iota^\dagger \msG \cong \iota^\dagger \iota_{!*}\msF \cong \msF$,
it follows  from  Proposition \ref{Pr1} that  $\msF$ is of type $\mcA$.

Finally, we shall consider (iii).
Suppose that $\msF$ is of type $\mcA$.
Let us consider the short exact sequence
\begin{align} \label{Eq222}
0 \longmigi \msF_1 \longmigi \msF \longmigi \msF/\msF_1 \longmigi 0.
\end{align}
All the $\mcD_X$-modules in this sequence are holonomic, so
there exists a dense open subscheme $U$ of $X$ such that their restrictions to $U$  are flat bundles.
Denote by $\eta$ the natural open immersion $U \migiincl X$ (hence $\eta^\star (-) \cong \eta^\dagger (-)$ on holonomic $\mcD_X$-modules,  see ~\cite[Theorem 2.7.1, (ii)]{HTT}).
Since $\msF/\msF_1$ is irreducible, the  morphism 
$\msF/\msF_1 \migi \int_\eta \eta^\dagger (\msF/\msF_1)$ induced by the adjunction relation ``$\eta^\star (-)\left(\cong \eta^\dagger (-) \right) \dashv \int_\eta (-)$"
is verified to be injective. 
Thus, by putting
$\msG$ as  the cokernel of this injection,  we obtain  a short exact sequence
\begin{align} \label{Eoa23}
0 \longmigi \msF/\msF_1 \longmigi \int_\eta \eta^\dagger (\msF/\msF_1) \longmigi \msG \longmigi 0.
\end{align}
Since $\eta^\dagger (\msF/\msF_1)  = \eta^\dagger \left( \int_\eta \eta^\dagger (\msF/\msF_1)\right)$,
the sequence (\ref{Eoa23}) implies 
 the equality $\eta^\dagger \msG = 0$,  meaning  that $\msG$ is supported in $X \setminus U$.
By 
Corollary \ref{Co22fe},
$\msG$ turns out to be of type $\mcA$.
On the other hand,  the  flat bundle    $\eta^\dagger \msF$  belongs to $\mcA (U)$,
 so the property ($\alpha$)  on $\mcA$  implies that
 its quotient flat bundle $\eta^\dagger (\msF/\msF_1)$ belongs to $\mcA (U)$.
 Hence, 
the direct image  $\int_\eta \eta^\dagger (\msF/\msF_1)$ is of type $\mcA$ (cf.  Proposition \ref{L011}, (ii)).
  By the exactness of  (\ref{Eoa23}) and  Proposition \ref{P01},  $\msF/\msF_1$ is verified to be    of type $\mcA$.
Moreover, it follows from  Proposition \ref{P01} again  and the exactness of (\ref{Eq222}) that
  $\msF_1 \in D_h^b (\mcD_X)$.
By applying successively an argument similar  to this, we see that every subquotient $\msF_j/\msF_{j+1}$ ($j=0, \cdots, m$) is of type $\mcA$.
This completes the proof of assertion (iii). 
\end{proof}
\SSP

\vspace{10mm}
\section{Holonomic $\mcD$-modules of arithmetic types} \LSP

In the rest of the present paper, 
we focus  on certain specific types of ``$\mcA$"'s, i.e., the fibered categories  of $G$-connections, globally nilpotent connections, and  almost everywhere nilpotent connections, respectively. 
We first recall their definitions and then
  introduce 
the  derived categories of $\mcD$-modules of such types.
The goal of this section is to prove Theorem \ref{Theer}.

For each commutative ring $R_0$, we shall write $\mbA^1_{R_0}$ for the affine line over $R_0$, i.e., $\mbA^1_{R_0} := \mr{Spec}(R_0 [t])$.

\LSP
\subsection{Global inverse radius of a flat bundle} \label{SS01}

Let $K$ be a number field, and 
denote by  $\mcO_K$ its ring of integers.
Consider a nonempty open subscheme $\mr{Spec}(R)$
  of $\mr{Spec}(\mcO_K)$ (where $R$ denotes a Dedekind domain with $\mcO_K \subseteq R \subseteq K$).
Denote by $\Sigma_R$ the set of  closed points of $\mr{Spec}(R)$, i.e., the set of finite places of $K$ having center on $R$.
For each
prime number $p$ and each
 $v \in \Sigma_R$ with $v |p$, we denote by $|-|_v$ the non-archimedean absolute value of $K$ corresponding to  
$v$, normalized as $|p|_v = p^{- [\widehat{K}_v: \mbQ_p]/[K: \mbQ]}$,
where $\widehat{K}_v$ denotes  the $v$-completion of $K$. 
Also, denote by $\widehat{\mcO}_v$ the ring of integers of $\widehat{K}_v$, and by $k (v)$ the residue field of $\widehat{K}_v$ of characteristic $p =: p (v)>0$.

Next, let $X_K$ be a smooth algebraic variety over $K$ 
and
 $f : X_R \migi \mr{Spec}(R)$ a smooth $R$-scheme  of finite type with geometrically connected non-empty fibers equipped with an isomorphism $X_R  \times_R K \isom X_K$;
 such an  $R$-scheme $X_R$ will be called a {\it model} of $X_K$ over $R$.
Denote by $K_{X_K}$ the function field of $X_K$ and by $d$ the relative dimension of $X_R/R$.
For each $v \in \Sigma_R$, there exists a unique extension $| - |_{X_R, v}$ of $| -|_v$ to a non-archimedean absolute value of  $K_{X_K}$
such that  the local ring $\mcO_{X_R, \eta_v}$ corresponding to
the generic point $\eta_v$ of the closed fiber $X_{k (v)} := X_R \times_R k(v)$
coincides with $\left\{ x \in K_{X_K} \, | \, |x|_{X_R, v} \leq 1 \right\}$.
For example, when $X_K = \mbA^1_K$, the absolute value on $K (t) \left(= K_{\mbA^1_{K}} \right)$ associated to both  the model $\mbA_{\mcO_K}^1$ of $\mbA^1_K$ and  $v \in \Sigma_{\mcO_K}$ coincides with 
 the Gauss absolute value, i.e., the absolute value   
  given by 
\begin{align} \label{dd00}
\left| \frac{\sum_i a_i t^i}{\sum_i b_i t^i}\right|_{\mr{Gauss}, v} := \frac{\sup_i |a_i|_v}{\sup_i |b_i|_v}.
\end{align} 
Also,   
$| - |_{X_R, v}$
naturally induces a norm $|\!| - |\!|_{X_R, v}$ on $M_{n \times n} (K_{X_K})$ (= the $K_{X_K}$-vector space    of $n \times n$ matrices with entries in $K_{X_K}$) given by $|\!| G |\!|_{X_R, v} := \mr{max}\left\{ | g_{i, j} |_{X_R, v} \, | \, 1 \leq i \leq n, 1 \leq j \leq n\right\}$ for any $G := (g_{i, j})_{i, j} \in M_{n \times n} (K_{X_K})$.

Let $\msF := (\mcF, \nabla)$ be a rank $n$ flat bundle defined on 
a nonempty open  subscheme of $X_K$.
We shall choose an \'{e}tale local (relative) coordinate $\underline{x}:= (x_1, \cdots, x_d)$ of  $X_R/R$ around
the point $\eta_v$ and a  basis $\underline{e} := (e_1, \cdots, e_n)$  of 
 $\mcF$ over that point.
 For any $\underline{\alpha} := (\alpha_1, \cdots, \alpha_d) \in \mbZ_{\geq 0}^d$,
we shall set
\begin{align} \label{Eq1}
\nabla_{[\underline{\alpha}]} := \prod_{i=1}^d  \nabla \left(\frac{\partial}{\partial x_i}\right)^{\alpha_i}.
\end{align}
Then, there exists  a unique  $n \times n$ matrix 
\begin{align} \label{Eq341}
G (\nabla)_{[\underline{\alpha}]} \in M_{n \times n} (K_X)
\end{align}
with entries in $K_X$  satisfying 
$\nabla_{[\underline{\alpha}]}
 \underline{e} = \underline{e} G (\nabla)_{[\underline{\alpha}]}$,
where  $G (\nabla)_{[\underline{0}]} := I_n$, i.e., the identity matrix of size $n$.

Recall that the {\bf radius of convergence  of $\msF$ at $v \in \Sigma_R$} is defined as the value
\begin{align} \label{Eq34}
\mr{Rad}_{X_R, v} (\msF) := \left(\mr{max} \left\{1, \limsup_{|\underline{\alpha}| \to \infty} \left|\!\left| \frac{1}{\underline{\alpha}!} \cdot G (\nabla)_{[\underline{\alpha}]} \right|\!\right|_{X_R, v}^{\frac{1}{|\underline{\alpha}|}} \right\} \right)^{-1} \in (0, 1],
\end{align}
 where $\underline{\alpha}! := \prod_{i=1}^d \alpha_i !$ and $|\underline{\alpha}| := \sum_{i=1}^d \alpha_i$.
 This value depends neither  on the choices of $\underline{x}$ nor $\underline{e}$ (cf. ~\cite[Proposition 2.9]{DVi} or ~\cite[Proposition 1.3]{ChDw}).
 Moreover, 
 by using these values for various elements $v \in \Sigma_R$, one may define the {\bf global inverse radius of $\msF$}   as 
\begin{align} \label{Eq33}
\rho_{X_R} (\msF) := \sum_{v \in \Sigma_R} \log \frac{1}{\mr{Rad}_{X_R, v} (\msF)} \in \mbR_{\geq 0} \sqcup \{ \infty \}.
\end{align}

\SSP
\begin{rema} \label{Rem4}
We here describe several  properties of radii  defined above.
\begin{itemize}
\item[(i)]
If $\msG$ is  another  flat bundle,
then the tensor product $\msF \otimes \msG$ of flat bundles $\msF$, $\msG$   satisfies 
\begin{align} \label{Eq104}
\mr{Rad}_{X_R, v} (\msF \otimes \msG) \geq \mr{min}\left\{ \mr{Rad}_{X_R, v} (\msF), \mr{Rad}_{X_R, v} (\msG) \right\}
\end{align}
 for every $v \in \Sigma_R$ (cf. ~\cite[Lemma 6.2.8, (c)]{Ked}).
This implies the inequality 
\begin{align} \label{Ed24}
\rho_{X_R} (\msF \otimes \msG) \leq \rho_{X_R} (\msF) +  \rho_{X_R} (\msG).
\end{align}
\item[(ii)]
Let $0 \migi \msE \migi  \msF \migi \msG \migi 0$ be a short exact sequence of flat bundles.
Then, it follows from ~\cite[Chap.\,IV, \S\,2.5, Proposition 1]{And1} that, for each $v \in \Sigma_R$,  the following inequality holds:
\begin{align} \label{Eq101}
\mr{Rad}_{X_R, v} (\msF) =  \mr{min} \left\{ \mr{Rad}_{X_R, v} (\msE), \mr{Rad}_{X_R, v} (\msG) \right\}.
\end{align}
Hence,   it is verified that
\begin{align} \label{Eq100}
\mr{max} \left\{\rho_{X_R}(\msE), \rho_{X_R}(\msG) \right\} \leq \rho_{X_R} (\msF) \leq \rho_{X_R}(\msE) + \rho_{X_R} (\msG).
\end{align}
\item[(iii)]
Let  $g : Y_R \migi \mr{Spec}(R)$ be another  smooth $R$-scheme of finite type with geometrically connected non-empty fibers.
Also, let $h_R : Y_R \migi X_R$ be a smooth $R$-morphism.
The pull-back of $\msF$ via  the generic fiber $h$ of $h_R$  specifies a flat bundle $h^*\msF$  on   an open subscheme of $Y_K := Y \times_R K$.
Then, we can prove the inequality
 \begin{align} \label{Eq102}
 \mr{Rad}_{Y_R, v} (h^*\msF) \geq \mr{Rad}_{X_R, v} (\msF)
 \end{align}
  for every $v \in \Sigma_R$.
This implies 
\begin{align} \label{Ed23}
\rho_{Y_R} (h^*\msF)\leq \rho_{X_R} (\msF).
\end{align}
\end{itemize}
\end{rema}

\LSP
\subsection{Arithmetic properties  on flat bundles} \label{SS0d1}


\SSP
\bde \label{Deg44}
Let $\msF := (\mcF, \nabla)$ be a flat bundle on $X_K/K$.
\begin{itemize}
\item[(i)]
We shall say that
$\msF$, or $\nabla$,  is  {\bf globally convergent}  (resp., {\bf of type $G$})
 if $\mr{Rad}_{X_R, v}(\msF) =1$ for all but finitely many elements $v \in \Sigma_R$  (resp., $\rho_{X_R} (\msF) < \infty$).
 \item[(ii)]
We shall say that $\msF$, or $\nabla$, is {\bf almost everywhere (a.e.)\,nilpotent} (resp., {\bf globally nilpotent}) if  
there exists a flat bundle $\msF_R$ on 
a nonempty open subscheme of $X_R$ relative to $\mr{Spec}(R)$ satisfying
the two conditions (1) and (2) (resp., (1) and (3)) described below:
\begin{itemize}
\item[(1)]
$\msF_R$ is isomorphic to $\msF$ at   the generic point of $X_K$;
\item[(2)]
Let $\Sigma_{\msF_R}^{\mr{nilp}}$ denote the set  of prime numbers $p$  such that the flat bundle $\msF_R \otimes_R k(v)$ induced by $\msF_R$  has nilpotent $p$-curvature for any $v \in \Sigma_R$, $v| p$.
 Then, $\Sigma_{\msF_R}^{\mr{nilp}}$  has Dirichlet density one (cf., e.g., ~\cite[\S\,5]{Kat} for the definition of $p$-curvature).
\item[(3)]
The flat bundles $\msF_R \otimes_R k(v)$  have nilpotent $p$-curvature for all but finitely many elements  $v \in \Sigma_R$.
\end{itemize}
\end{itemize}
\ede
\SSP

\begin{rema} \label{Effd}
It is immediate that the various definitions described  in Definition \ref{Deg44} are independent of the choices of the Dedekind domain $R$ and the model $X_R$ of $X_K$.
\end{rema}
\SSP

Next, we shall fix  a smooth algebraic variety $X$ over $\overline{\mbQ}$.

\SSP
\bde \label{Def1}
Let 
$\msF := (\mcF, \nabla)$   be 
a flat bundle on $X/\overline{\mbQ}$.
We shall say that $\msF$ is 
{\bf globally convergent} (resp., {\bf of type $G$} ; resp., {\bf a.e.\,nilpotent}; resp., {\bf globally nilpotent})
if there exists a smooth  algebraic variety $X_K$ over a number field $K$
and 
a flat bundle $\msF_K$ on $X_K/K$ such that
$(X_K, \msF_K) \times_K \overline{\mbQ} \cong (X, \msF)$ and that $\msF_K$ is globally convergent  (resp., of type $G$; resp., a.e.\,nilpotent; resp., globally nilpotent) in the sense of Definition \ref{Deg44}.
 If $\msF$ is of type $G$, then we sometimes call $\nabla$ a {\bf $G$-connection}.
\ede
\SSP

Regarding the various   notions mentioned above, we have the following implications for an arbitrary flat bundle $\msF$ (as mentioned in ~\cite[\S\,1.4]{AB}):
\begin{align} \label{Eqww98}
\text{$\msF$ is globally convergent}
 & \Longrightarrow
\text{$\msF$ is of type $G$ (resp., globally nilpotent)} \\
 &\Longrightarrow
\text{$\msF$ is a.e.\,nilpotent}. \notag
\end{align}

 The following assertion is a direct consequence  of the previous study  by T. Honda concerning   the Grothendieck-Katz $p$-curvature for  the rank one case.
 
\SSP
\bpr \label{Prop562}
Suppose that $X = X_R \times_R \overline{\mbQ}$ for some 
open subscheme $X_R$ of  $\mbA_R^1$  (where $R$ is as above).
Also, let $\msF_R$ be a rank one  flat bundle on $X_R/R$, and write  $\msF := \msF_R \times_R \overline{\mbQ}$.
Then, the conditions described in (\ref{Eqww98}) are equivalent to each other, and moreover, 
these are equivalent to each of the following two conditions:
\begin{itemize}
\item[($*$)]
 The 
 flat bundles $\msF_R \otimes_R k (v)$ 
 have vanishing $p (v)$-curvature  for
all but finitely many
  elements  $v$ of  $\Sigma_R$.
\item[($**$)]
$\msF$ becomes generically  trivial (i.e., trivial at the generic point) after pulling back via 
  a finite \'{e}tale  covering of $X$.
\end{itemize}
\epr
\begin{proof}
Suppose that
$\msF$ is a.e.\,nilpotent.
Then,
the set $\Sigma_{\msF_R}^{\mr{nilp}}$ (cf. Definition \ref{Deg44}, (ii)) has Dirichlet density one.
Since $\msF_R$ is of rank one,  
the flat bundle $\msF_R \otimes_R k (v)$ has vanishing $p (v)$-curvature for every $v$ with $p (v) \in \Sigma_{\msF_R}^{\mr{nilp}}$.
Here, we recall the proof of the Grothendieck-Katz conjecture for rank one flat bundles on $X$ given by T. Honda (cf. ~\cite{Hon}), in which he proved that  the assertion of this conjecture is equivalent to Chebotarev's density theorem;
in particular, by applying that discussion to our situation, we see that the condition ($**$) is satisfied.
The implication  ($**$) $\Rightarrow$ ($*$) is clear and well-known.
Finally,  it follows from 
~\cite[Proposition 3.3]{DVi}   that ($*$) implies the global convergency of $\msF$.
This completes the proof of this proposition.
\end{proof}

\LSP
\subsection{Holonomic $\mcD$-modules  of arithmetic types} \label{SS0ww1}

We shall set
\begin{align}
\mcA_{\overline{\mbQ}, G} \ \left(\text{resp.,} \  \mcA_{\overline{\mbQ}, \mr{nilp}}; \text{resp.},  \  \mcA_{\overline{\mbQ}, \mr{aen}} \right)
\end{align}
to be the fibered subcategory  of $\mcF \mcB_{\overline{\mbQ}}$
classifying 
flat bundles of type $G$ (resp., globally nilpotent flat bundles; resp.,  a.e.\,nilpotent   flat bundles).
The following assertion in the case of $\mcA_{\overline{\mbQ}, G}$ is a direct consequence of  the main result of  ~\cite{AB}.

\SSP
\bpr \label{Thr4}
The fibered subcategories  $\mcA_{\overline{\mbQ}, G}$, 
$\mcA_{\overline{\mbQ}, \mr{nilp}}$, and  $\mcA_{\overline{\mbQ}, \mr{aen}}$
 satisfy  the three properties  ($\alpha$)-($\gamma$) described in  \S\,\ref{SS0233}.
\epr
\begin{proof}
Let $\mcA \in \left\{\mcA_{\overline{\mbQ}, G},  \mcA_{\overline{\mbQ}, \mr{nilp}}, \mcA_{\overline{\mbQ}, \mr{aen}} \right\}$.
As mentioned in ~\cite[\S\,1.4]{AB},
the  subcategory $\mcA(Y)$  of $\mcF \mcB_{\overline{\mbQ}}(Y)$ determined by  each $Y \in \mr{ob}(\mcS m_{\overline{\mbQ}})$
is closed  under various operations so that $\mcA$ satisfies ($\alpha$).
Also,  the property   ($\gamma$) can be verified from the definitions involved.
Finally, the property  ($\beta$) follows from ~\cite[Main Theorem]{AB} (if $\mcA = \mcA_{\overline{\mbQ}, G}$) and ~\cite[Theorem 5.10]{Kat} (if $\mcA =   \mcA_{\overline{\mbQ}, \mr{nilp}}$ or $\mcA_{\overline{\mbQ}, \mr{aen}}$).
\end{proof}
\SSP

By applying Proposition \ref{Thr4} and the arguments in the previous section,
we obtain  the triangulated subcategory 
\begin{align} \label{Eq2201}
D^b_h (\mcD_X)_G := D^b_h (\mcD_X)_{\mcA_{\overline{\mbQ}, G}}\hspace{40mm} 
 \\
  \left(\text{resp.,} \ 
  D^b_h (\mcD_X)_{\mr{nilp}} :=D^b_h (\mcD_X)_{\mcA_{\overline{\mbQ}, \mr{nilp}}} ; \text{resp.,} \   D^b_h (\mcD_X)_{\mr{aen}} :=D^b_h (\mcD_X)_{\mcA_{\overline{\mbQ}, \mr{aen}}}  \notag 
  \right)
\end{align}
of  $D_h^b (\mcD_X)$.
The implications between the conditions displayed in (\ref{Eqww98}) give rise to  inclusion relations of categories
\begin{align} \label{Eqkkid}
D^b_h (\mcD_X)_G \subseteq D^b_h (\mcD_X)_{\mr{aen}},  \ \ \  
D^b_h (\mcD_X)_{\mr{nilp}} \subseteq D^b_h (\mcD_X)_{\mr{aen}}.
\end{align}
In particular, by Theorem \ref{Th40} applied to $\mcA \in \left\{ \mcA_{\overline{\mbQ}, G}, \mcA_{\overline{\mbQ}, \mr{nilp}}, \mcA_{\overline{\mbQ}, \mr{aen}}\right\}$,
we obtain six functors on the respective  derived categories.
This proves Theorem \ref{Theer}.

\SSP
\begin{rema} \label{Rem7}
Let $\Box \in \{ G,  \mr{nilp}, \mr{aen}\}$, and 
denote by $D_{rh}^b (\mcD_X)$ the full subcategory of $D_h^b (\mcD_X)$ consisting of chain complexes having {\it regular} holonomic cohomology.
If $\msF$ is an a.e.\,nilpotent flat bundle on $X/\overline{\mbQ}$,
then  it follows from a result by Katz (cf. ~\cite[Theorem 8.1]{Kat4} or ~\cite[Theorem 6.1 and Remark 6.3]{DGS}) that
it has at most regular singularities at infinity along any smooth curve in $X$ and has rational exponents.
This fact together with (\ref{Eqww98}) implies that, for  any  chain complex $\msF^\bullet$   in $D_h^b(\mcD_X)_\Box$,  
the cohomology sheaf $\mcH^l (\msF^\bullet)$ (for each $l \in \mbZ$) defines  a regular holonomic $\mcD_X$-module.
In particular,  the  inclusion 
$D_h^b (\mcD_X)_\Box \subseteq D_h^b (\mcD_X)$
  factors through $D_{rh}^b (\mcD_X) \subseteq D_h^b (\mcD_X)$, i.e., we have
  \begin{align} \label{Eq64}
  D_h^b (\mcD_X)_\Box \subseteq D_{rh}^b (\mcD_X).
  \end{align}
\end{rema}

\LSP
\subsection{Global inverse radius of a holonomic $\mcD_{\mbA^1}$-modules} \label{SS0144}


We have extended the class of $G$-connections to chain complexes in $D^b_h (\mcD_{\mbA^1})$.
The global inverse radius  can be generalized to
 an invariant  $\rho (\msF^\bullet) \in \mbR_{\geq 0} \sqcup \{ \infty \}$ associated to each $\msF^\bullet \in \mr{ob}(D^b_h (\mcD_{(-)})_G)$ accordingly.
 For simplicity,  we only consider the case where $X$ is the affine line $\mbA^1 := \mbA^1_{\overline{\mbQ}}$ over $\overline{\mbQ}$.
 
First, let $\msF$ be a holonomic $\mcD_{\mbA^1}$-module   that belongs 
 to $D^b_h (\mcD_{\mbA^1})_G$  as a complex concentrated at degree $0$.
There exist
a number field $K$ and
a dense open subscheme $U$ of $\mbA^1_{\mcO_K}$
such that the restriction $\msF |_{U \times_{\mcO_K} \overline{\mbQ}}$ of $\msF$ to the open subscheme $U \times_{\mcO_K} \overline{\mbQ} \left( \subseteq \mbA^1\right)$ can be  obtained as the pull-back of a flat bundle $\msF_{K}$ on $(U \times_{\mcO_K}K)/K$.
 Then, the  value $\rho (\msF):= \rho_{\mbA^1_{\mcO_K}} (\msF_K)$ (cf. (\ref{Eq33})) can be defined and  depends neither  on the choices of $K$, $U$, nor $\msF_K$ (cf. ~\cite[\S\,1.3]{AB}).
 More generally, for each $\msF^\bullet \in D^b_h (\mcD_{\mbA^1})_G$,
  we  set 
 \begin{align} \label{dd32}
 \rho (\msF^\bullet) := \mr{max} \left\{ \rho (\mcH^l (\msF^\bullet)) \, | \, l \in \mbZ \right\}.
 \end{align}
In this way, we have obtained a well-defined map 
\begin{align}
\rho : \mr{ob}(D^b_h (\mcD_{\mbA^1})_G) \migi \mbR_{\geq 0}.
\end{align}
extending $\rho_{\mbA^1_{(-)}}$.

We shall generalize the second inequality in (\ref{Eq100}), as follows.

\SSP
\bpr \label{Pr34}
Let  $\msE^\bullet \migi \msF^\bullet \migi \msG^\bullet \xrightarrow{+1}$ be  a distinguished triangle of chain complexes in $\mcD_h^b (\mcD_{\mbA^1})_G$.
 Then, the following inequality holds:
\begin{align} \label{Eq92}
\left| \rho (\msE^\bullet) - \rho (\msG^\bullet) \right| \leq \rho (\msF^\bullet) \leq \rho (\msE^\bullet) + \rho (\msG^\bullet).
\end{align}
\epr
\begin{proof}
The distinguished triangle $\msE^\bullet \migi \msF^\bullet \migi \msG^\bullet \xrightarrow{+1}$ induces a long  exact sequence of $\mcD_{\mbA^1}$-modules of type $G$:
\begin{align} \label{Ekkid2}
\cdots \xrightarrow{\beta_{l-1}}
\mcH^{l-1}(\msG^\bullet) \xrightarrow{\gamma_{l-1}}
\mcH^{l}(\msE^\bullet)\xrightarrow{\alpha_l}
\mcH^{l}(\msF^\bullet)\xrightarrow{\beta_l}
\mcH^{l}(\msG^\bullet)\xrightarrow{\gamma_l}
\mcH^{l+1}(\msE^\bullet)\xrightarrow{\alpha_{l+1}} \cdots.
\end{align}
To prove  the second inequality of  (\ref{Eq92}), let us consider the following short exact sequences arising from  (\ref{Ekkid2}):
\begin{align}
 \label{Eq105}0 \migi \mr{Coker}(\gamma_{l-1}) \left(= \mr{Im}(\alpha_l) \right)
&\migi \mcH^l (\msF^\bullet) \migi \mr{Im}(\beta_l) \left(=  \mr{Coker}(\alpha_l)\right) \migi 0, \\
\label{Eq106}0 \migi  \mr{Im}(\gamma_{l-1}) &\migi \mcH^{l}(\msE^\bullet) \migi  \mr{Coker}(\gamma_{l-1})  \migi 0, \\
\label{Eq107}
0 \migi \mr{Im}(\beta_l) &\migi \mcH^{l}(\msG^\bullet) \migi \mr{Coker}(\beta_l ) \migi 0.
\end{align}
By the fact mentioned in Remark \ref{Rem4}, (ii), 
these 
   induce    inequalities 
\begin{align}
\rho (\mcH^l (\msF^\bullet)) &\leq \rho (\mr{Coker}(\gamma_{l-1})) + \rho (\mr{Im}(\beta_l)), \\
\rho (\mr{Coker}(\gamma_{l-1})) &\leq \rho (\mcH^l (\msE^\bullet)), \text{and}\notag \\
\rho (\mr{Im}(\beta_l)) &\leq \rho (\mcH^l (\msG^\bullet)), \notag 
\end{align}
respectively.
This  implies $\rho (\mcH^l (\msF^\bullet)) \leq \rho (\mcH^l (\msE^\bullet)) + \rho (\mcH^l (\msG^\bullet))$, and hence, we obtain  the second inequality of  (\ref{Eq92}), as desired.

Also, entirely similar arguments enable us to obtain  inequalities 
$\rho (\msE^\bullet) \leq \rho (\msG^\bullet) + \rho (\msF^\bullet)$ and $\rho (\msG^\bullet) \leq \rho (\msF^\bullet) + \rho (\msE^\bullet)$,  which together are equivalent to the first inequality of (\ref{Eq92}).
This completes the proof of this assertion.
\end{proof}

\vspace{10mm}
\section{Middle convolution on holonomic $\mcD$-modules of arithmetic types}  \label{S1132}
\LSP

This section discusses  the middle convolution functors (with rational parameters) on 
the derived categories under consideration.
As a consequence, we prove Theorem  \ref{ThCC}.

\LSP
\subsection{Middle convolution} \label{SS023}

Let $q$ be a $\overline{\mbQ}$-rational point  of  $\mbA^1 := \mbA^1_{\overline{\mbQ}} \left(= \mr{Spec}(\overline{\mbQ}[t])\right)$;
it determines an open immersion $\eta^q : \mbA^1 \setminus \{q \} \migiincl \mbA^1$.
We will use the same notation ``$q$" to denote the corresponding element in  $\overline{\mbQ}$.
For  
each 
 $\lambda \in \overline{\mbQ}$,
 we shall set
\begin{align} \label{Ewrrt}
\msK^\lambda_q := \eta^q_{*}\left(\mcO_{\mbA^1 \setminus \{ q\}}, d + \lambda \cdot \frac{d(t-q)}{t-q}\right),
\end{align}
which is a holonomic $\mcD_{\mbA^1}$-module; it is irreducible   when $\lambda \notin \mbZ$.
Up to isomorphism, $\msK^\lambda_q$ depends only on the image of $\lambda$ in $\overline{\mbQ}/\mbZ$.

\SSP
\bpr \label{Prop443}
Let $\Box \in \{ G, \mr{nilp}, \mr{aen} \}$.
Then, the  $\mcD_{\mbA^1}$-module $\msK^\lambda_q$ belongs to  $D^b_h (\mcD_{\mbA^1})_\Box$ if and only if $\lambda \in \mbQ$.
\epr
\begin{proof}
Let $K$ and $R$ be as in \S\,\ref{SS01}; we may assume (after possibly replacing these with different ones) that $\lambda, q \in R$.
Write $\nabla_R$ for the $R$-connection $d + \lambda \cdot \frac{d(t-q)}{t-q}$ on
$\mcO_{X_R}$, where
 $X_R := \mbA^1_R \setminus \overline{\{ q \}}$, and write $(X, \nabla) := (X_R, \nabla_R) \times_R \overline{\mbQ}$.

Now, let us consider the ``if" part of the required equivalence.
Suppose that $\lambda \in \mbQ$.
For a prime number  $p$ and an element $v \in \Sigma_R$
with $v | p$,
one  can reduce  $\nabla_R$ to obtain a $k(v)$-connection $\nabla_{R, v}$ on $\mcO_{X_R \otimes_R k (v)}$.
The reduction of $\lambda$ modulo $p$ belongs to $\mbF_p := \mbZ/p\mbZ$.
Hence, if $C$ denotes the Cartier operator on $\Omega_{X_R \otimes_R k (v)/k (v)}$, then the equality $C (\lambda \cdot \frac{d (t-q)}{t-q}) = \lambda \cdot \frac{d (t-q)}{t-q}$ holds.
According to ~\cite[Corollary 7.1.3]{Kat3},
 the connection  $\nabla_{R, v}$  has vanishing $p$-curvature, so 
  Proposition \ref{Prop562}
 implies  that
 $(\mcO_{X}, \nabla)$ is    globally convergent.
By the implications in  (\ref{Eqww98}), 
this  flat bundle is also of type $G$ (resp., globally nilpotent; resp., a.e.\,nilpotent).
 It follows from Proposition \ref{L011}, (i),   that $\msK^\lambda_q \left(\cong \int_{\eta^q} (\mcO_{X }, \nabla) \right)$ is  contained in  $D_h^b (\mcD_{\mbA^1})_G$ (resp., $D_h^b (\mcD_{\mbA^1})_{\mr{nilp}}$; resp., $D_h^b (\mcD_{\mbA^1})_{\mr{aen}}$).
This completes the proof of the ``if" part.

Next, we shall  prove 
 the inverse direction of the required equivalence.
 To this end,
 it suffices to consider the case where ``$\Box = \mr{aen}$"  because of the implications in  (\ref{Eqww98}).
Suppose that $\msK^\lambda_q$ belongs to $D^b_h (\mcD_{\mbA^1})_{\mr{aen}}$.
 Since $\eta^{q\dagger}(\msK^\lambda_q) \cong (\mcO_{X}, \nabla)$,
 it follows from  Proposition \ref{Pr1} that $(\mcO_{X}, \nabla)$ is a.e.\,nilpotent.
 Hence, by Proposition \ref{Prop562},
 $(\mcO_{X}, \nabla)$ is verified to be globally covergent.
 By reversing the steps in  the above discussion,
 we see that 
 the mod $v$ reduction of $\lambda$ lies in $\mbF_{p (v)}$ for every element $v$ in $\Sigma_R$.
 This implies
 $\lambda \in \mbQ$ by Chebotarev's density theorem, thus completing 
    the proof of the ``only if" part.
\end{proof}
\SSP

Denote by $\mbP^1 \left(\supseteq \mbA^1 \right)$ the projective line over $\overline{\mbQ}$ and by $\pi_1$ (resp., $\pi_2$)  the first projection $\mbA^1 \times_{\overline{\mbQ}} \mbA^1 \migi \mbA^1$ (resp., the second projection $\mbP^1 \times_{\overline{\mbQ}} \mbA^1 \migi \mbA^1$).
Also, denote by $\mu : \mbA^1 \times_{\overline{\mbQ}} \mbA^1 \migi \mbA^1$  the morphism given by $(x, y) \mapsto y-x$ and  by $\eta_{\mbA} : \mbA^1 \times_{\overline{\mbQ}} \mbA^1 \migiincl \mbP^1 \times_{\overline{\mbQ}} \mbA^1$ the natural open immersion.
For a
chain complex  $\msF^\bullet$ in  $D^b_h (\mcD_{\mbA^1})$
 and $\lambda \in \overline{\mbQ} \setminus \mbZ$,
 we shall set
\begin{align} \label{Eq400k}
\mr{mc}_\lambda (\msF^\bullet) := \int_{\pi_2}\eta_{\mbA !*}(L\pi_1^* \msF \otimes^L_{\mcO_{\mbA^1 \times \mbA^1}} L\mu^* \msK_0^\lambda) \in D_h^b(\mcD_{\mbA^1})
\end{align} 
(cf. ~\cite[\S\,2.8]{Kat2} for the corresponding definition in the $l$-adic setting).
 Also, for each $l \in \mbZ$, we shall set
\begin{align} \label{Eq82}
\mr{mc}_\lambda^l (\msF^\bullet) :=  \mcH^l (\mr{mc}_\lambda (\msF^\bullet)).
\end{align}
 We  refer to $\mr{mc}_\lambda (\msF^\bullet)$ as the {\bf (additive) middle convolution of $\msF$ with the parameter $\lambda$}.
 The resulting assignment $\msF^\bullet \mapsto \mr{mc}_\lambda (\msF^\bullet)$ defines an endofunctor   on $D_h^b (\mcD_{\mbA^1})$:
 \begin{align}  \label{Eq81}
 \mr{mc}_\lambda : D_h^b (\mcD_{\mbA^1}) \migi  D_h^b (\mcD_{\mbA^1}).
 \end{align}

In the case of flat bundles, the definition of middle convolution just mentioned coincides with the definition given in ~\cite{Ari} (cf. ~\cite[Lemma 6.9]{Ari2}); see  the former assertion of Theorem \ref{TT5} described later.

The following assertion is a direct consequence of results obtained so far (cf.  ~\cite[Theorem 1]{DR2} for the case of globally nilpotent connections).

\SSP
\bt \label{T491}
Let $\Box \in \{ G, \mr{nilp}, \mr{aen} \}$, and 
 suppose that $\lambda$ belongs to  $\mbQ \setminus \mbZ$.
Then,
the functor $\mr{mc}_\lambda$  restricts to an endofunctor
\begin{align} \label{Eqwo892}
\mr{mc}_\lambda : D_h^b (\mcD_{\mbA^1})_\Box \migi  D^b_h (\mcD_{\mbA^1})_\Box
\end{align}
on $D_h^b (\mcD_{\mbA^1})_\Box$.
\et
\begin{proof}
 It follows from Proposition \ref{Prop443} that $\msK^\lambda_0$  is contained in  $D_h^b (\mcD_{\mbA^1})_\Box$.
Hence, by the definition of $\mr{mc}_\lambda (-)$  together with
Theorem \ref{Th40},  Propositions \ref{Prop909}, (i), and  \ref{Thr4},
we have  
$\mr{mc}_\lambda (\msF^\bullet) \in \mr{ob}(D_h^b (\mcD_{\mbA^1})_\Box)$ whenever $\msF^\bullet$ lies in $D_h^b (\mcD_{\mbA^1})_\Box$.
\end{proof}

\LSP
\subsection{Estimate of global inverse radii I} \label{SS023}

The rest of this section is devoted to 
  estimating the effect of  middle convolution on the global inverse radii.
 We first deal with  $\mcD_{\mbA^1}$-modules supported on a single point.
For each prime $p$ and $\lambda \in \overline{\mbQ} \setminus \mbZ$,
we shall set
\begin{align} \label{Ewiod}
\mr{ord}_p \lambda := - \sum_{v \in \Sigma_{\mcO_K}, v | p}\mr{log}_p |\lambda |_v,
\end{align}
where $K$ denotes a number field with $\lambda \in K$;
this value is immediately verified to be independent of the choice of $K$ and satisfies $\mr{ord}_p p =1$.
Moreover, we shall set
\begin{align} \label{df3}
H(\lambda) := \sum_{\substack{p :  \,\text{prime} \\ \text{s.t.\,$\mr{ord}_p \lambda < 0$}}} \left(\frac{1}{p-1} - \mr{ord}_p \lambda \right)\cdot \log p \ \left(>0 \right).
\end{align}


\SSP
\bpr \label{Prop558}
Denote by $\msO$ the trivial $\mcD_{\mr{Spec}(\overline{\mbQ})}$-module associated to the  $\overline{\mbQ}$-vector space $\overline{\mbQ}$.
\begin{itemize}
\item[(i)]
Let $\lambda$ be an element of $\overline{\mbQ}\setminus \mbZ$ and  
$q$ a $\overline{\mbQ}$-rational point  
of $\mbA^1$.
Then, we have
\begin{align}
\mr{mc}_\lambda (\int_q\msO) \cong \msK^\lambda_q.
\end{align}
\item[(ii)]
Let $\lambda$ be an element of  $\mbQ \setminus \mbZ$ and $(q_1, \cdots, q_n)  $ (where $n \in \mbZ_{>0}$) an $n$-tuple of $\overline{\mbQ}$-rational points of $\mbA^1$.
Then,  
the following inequality holds:
\begin{align}
\rho (\mr{mc}_\lambda (\bigoplus_{i=1}^n\int_{q_i}\msO)) \left(= \rho (\bigoplus_{i=1}^n\msK^\lambda_{q_i})\right) \leq H (\lambda).
\end{align}
\end{itemize}
\epr
\begin{proof}
First,  we shall consider assertion (i).
Denote by $q_\mbA$ the composite closed immersion $\mbA^1 \isom \{q \} \times_{\overline{\mbQ}}\mbA^1 \migiincl \mbA^1 \times_{\overline{\mbQ}} \mbA^1$.
The composite $\mu \circ q_\mbA$ coincides with the automorphism of $\mbA^1$ given by $x \mapsto x-q$.
This implies $L(\mu \circ q_\mbA)^* (\msK^\lambda_0) \cong
 (\mu \circ q_\mbA)^* (\msK^\lambda_0) \cong 
 \msK^\lambda_q$.
Hence, we have
\begin{align}
 \int_{\pi_2}\eta_{\mbA!*}(L\pi_1^*\int_{q} \msO\otimes^L L\mu^* \msK_0^\lambda)
 &\cong \int_{\pi_2}\eta_{\mbA!*}(\left(\int_{q_\mbA} (\msO_{\mbA^1}, d) \right)\otimes^L L\mu^* \msK^\lambda_0)   \\
 & \cong \int_{\pi_2} \eta_{\mbA!*}(\int_{q_\mbA}L(\mu \circ q_{\mbA})^* (\msK^\lambda_0))  \notag \\
 & \cong \int_{\pi_2} \eta_{\mbA!*}(\int_{q_\mbA}\msK^\lambda_q) \notag \\
 & \cong \msK^\lambda_q, \notag 
\end{align}
where 
\begin{itemize}
\item
the first ``$\cong$" follows from the base-change theorem (cf.  ~\cite[Theorem 1.7.3]{HTT}),
\item
 the second ``$\cong$" follows from the projection formula (cf. ~\cite[Corollary 1.7.5]{HTT}), and
 \item
  the last ``$\cong$" follows from some properties of minimal extensions asserted   in ~\cite[Theorem 3.4.2]{HTT} together with the simplicity of  $\msK^\lambda_q$.
\end{itemize}
This completes the proof of assertion (i).

Next, let us consider  assertion (ii).
Suppose that $\lambda \in \mbQ \setminus \mbZ$.
Also, let us take a number field $K$ with $q \in K$.
For each prime $p$ and each $v \in \Sigma_{\mcO_K}$ with $v | p$, 
we have 
\begin{align} \label{HHH4}
\limsup_{s \to \infty}\left| \frac{1}{s!} \cdot \prod_{j=0}^{s-1} (\lambda - j) \right|_{v}^{\frac{1}{s}} \leq 
 \begin{cases} 
 1 & \text{if $\mr{ord}_p \lambda \geq  0$}; \\
p^{\left(\frac{1}{p-1}-\mr{ord}_p \lambda\right) \cdot \frac{[\widehat{K}_v: \mbQ_p]}{[K: \mbQ]}}
& \text{if $\mr{ord}_p \lambda < 0$}
\end{cases}
\end{align}
 (cf. ~\cite[Theorem 47.8]{Sch}).
It follows that
\begin{align} \label{HHH5}
& \ \ \ \ \mr{Rad}_{\mbA^1_{\mcO_K}, v} (\bigoplus_{i=1}^n\msK_{q_i}^\lambda)^{-1} \\
&= \mr{max} \left(\left\{1\right\} \cup \bigcup_{i=1}^n \left\{ \limsup_{s \to \infty} \left|  \frac{1}{s!} \cdot \frac{1}{(t- q_i)^s} \cdot \prod_{j=0}^{s-1} (\lambda - j)\right|_{\mr{Gauss}, v}^{\frac{1}{s}} \right\}\right)
\notag  \\
&=  \mr{max}\left(\left\{1\right\} \cup \bigcup_{i=1}^r \left\{ \frac{1}{\mr{max} \left\{1, |q_i |_v \right\}} \cdot  \limsup_{s \to \infty} \left|  \frac{1}{s!} \cdot \prod_{j=0}^{s-1} (\lambda - j)\right|_{v}^{\frac{1}{s}} \right\}\right) \notag \\
&
\hspace{-2mm}
\stackrel{(\ref{HHH4})}{\leq} 
 \begin{cases} 
 1 & \text{if $\mr{ord}_p \lambda \geq  0$}; \\
p^{\left(\frac{1}{p-1}-\mr{ord}_p \lambda\right) \cdot \frac{[\widehat{K}_v: \mbQ_p]}{[K: \mbQ]}}
& \text{if $\mr{ord}_p \lambda < 0$}.
\end{cases}
\notag
\end{align}
Hence, we have 
\begin{align} 
\rho_{\mbA^1_{\mcO_K}} (\bigoplus_{i=1}^n \msK_{q_i}^\lambda)
&= \sum_{\text{$p$\,:\,prime}} \log \left(\prod_{v \in \Sigma_{\mcO_K}, v | p} \mr{Rad}_{\mbA^1_{\mcO_K}, v} (\bigoplus_{i=1}^n\msK_{q_i}^\lambda)^{-1}\right)   \\
& \hspace{-1.5mm} \stackrel{(\ref{HHH5})}{\leq} 
 \sum_{\substack{\text{$p$\,:\,prime} \\ \text{ s.t.\,$\mr{ord}_p\lambda < 0$}}}
  \log
\left(p^{\frac{1}{p-1}-\mr{ord}_p \lambda}\right)
\notag \\
& =  H (\lambda).\notag
\end{align}
This completes the proof of the latter assertion.
\end{proof}

\LSP
\subsection{Estimate  of global inverse radii II} \label{SS028v}

Next, let us consider  flat bundles on an open subscheme of $\mbA^1$.
As explained  in Remark \ref{REM99} described  later,  the general case can be investigated  by combining  Proposition \ref{Prop558} and the following theorem.

\SSP
\bt \label{TT5}
Let  $\lambda$ be an element of $\mbQ \setminus \mbZ$, $U$  a nonempty  open subscheme of $\mbA^1$,  and $\msF$ a  flat bundle on $U/\overline{\mbQ}$ of rank $n  \in \mbZ_{> 0}$ and type $G$.
Denote by $\iota$ the natural open immersion $U \migiincl \mbA^1$.
Then, both 
$\mr{mc}_\lambda (\int_\iota \msF) |_U$ and 
$\mr{mc}_\lambda (\iota_{!*} \msF) |_U$ are concentrated only at degree $0$ in $D_h^b (\mcD_{U})_G$ (hence these can be considered as $\mcD_U$-modules).
Moreover,  the following inequalities hold:
\begin{align} \label{Er33}
\left| \rho (\mr{mc}_\lambda (\int_\iota \msF) ) -\rho (\mr{mc}_\lambda (\iota_{!*}\msF)) \right| \leq  H (\lambda),  
\hspace{5mm} 
\rho (\mr{mc}_\lambda (\int_\iota \msF) ) \leq (n^2 +1) \left(\rho ( \msF) + H (\lambda)\right).
\end{align}

\et
\begin{proof}
To begin with, we shall prove the former assertion  and the first  inequality in (\ref{Er33}).
Denote by $\msG$ the cokernel of the natural inclusion  $\iota_{!*}\msF \migiincl \int_\iota \msF$.
Since $\msG$ is supported on $\mbA^1 \setminus U$,
it follows from Proposition \ref{Prop558}, (i), that
$\mr{mc}_\lambda^l (\msG) = 0$ for $l \neq 0$ and that 
$\mr{mc}_\lambda^0 (\msG)$ is isomorphic to a direct sum of finitely many $\msK^\lambda_q$'s (for some $q$'s).
Hence, Proposition \ref{Prop558}, (ii),  implies 
\begin{align} \label{E49}
\rho (\mr{mc}_\lambda (\msG)) \leq H (\lambda).
\end{align}
Also,  for each $x \in U (\overline{\mbQ})$ and $l \in \mbZ \setminus \{ 0\}$,  
the fiber of  $\mr{mc}_\lambda^l (\int_{\iota} \msF)$ over $x$ satisfies
\begin{align}
\mr{mc}_\lambda^l (\int_{\iota} \msF)_x \cong H^{l+1}_{\mr{dR}}(\mbP^1, \eta_{*}(\int_\iota \msF \otimes \msK^\lambda_x)) = 
H^{l+1}_{\mr{dR}}(\mbA^1, \int_\iota \msF \otimes \msK^\lambda_x)
=  0,
\end{align}
where
\begin{itemize}
\item
$\eta$ denotes the natural open immersion $\mbA^1\migiincl \mbP^1$ and the first ``$\cong$" follows  essentially from ~\cite[Corollary 2.8.5]{Kat2} and ~\cite[Theorem 7.1.1]{HTT} because 
$\msF$ has at most regular singularities
(cf. Remark \ref{Rem7}),
\item
the last ``$=$" for $l = -1$ follows from the fact that, since $x \in U (\overline{\mbQ})$, there is no horizontal section of $ \int_\iota \msF \otimes \msK^\lambda_x$ (on an open neighborhood of $x$), and 
\item
the last ``$=$" for $l \neq -1, 0$ follows from the affineness of   $\mbA^1$.
\end{itemize}
 
Similarly,
we have $\mr{mc}_\lambda^l (\iota_{!*}\msF)_x = 0$ for every $l \in \mbZ \setminus \{ 0\}$. 
It follows that $\mr{mc}_\lambda^l (\int_\iota \msF) |_U = \mr{mc}_\lambda^l (\iota_{!*} \msF) |_U = 0$ for every $l \in \mbZ \setminus \{ 0\}$.
In particular,  we have finished the proof of the former assertion.
Moreover, the natural short exact sequence $0 \migi \iota_{!*}\msF  \migi \int_\iota \msF \migi  \msG \migi 0$ induces 
  a short exact sequence of $\mcD_{U}$-modules
 \begin{align}
  0
  \longmigi   \mr{mc}_\lambda^0 (\iota_{!*}\msF) |_U
 \longmigi \mr{mc}_\lambda^0 (\int_{\iota} \msF)|_U
\longmigi 
 \mr{mc}_\lambda^0 (\msG)|_U
 \longmigi  0.
 \end{align}
 By  Proposition \ref{Pr34}
 and (\ref{E49}), we obtain  the inequalities
 \begin{align}
\left|  \rho (\mr{mc}_\lambda^0 (\int_\iota \msF) ) - \rho (\mr{mc}_\lambda^0 (\iota_{!*}\msF))  \right|
\leq 
  \rho (\mr{mc}_\lambda^0 (\msG)) \leq  H(\lambda).
 \end{align}
 This completes the proof of the first inequality in  (\ref{Er33}).
 
 Next, we shall prove the second inequality in  (\ref{Er33}).
 Let us take a number field $K$ such that 
 there exists an open subscheme $U_K$ of $\mbA^1_{K}$ and a flat bundle $\msF_K$ on $U_K/K$ with  $(U_K, \msF_K) \times_K \overline{\mbQ}   = (U, \msF)$.
  Denote by $\kappa : (U \times U) \setminus \Delta_U \migiincl \mbP^1 \times \mbA^1$ the natural open immersion, where $\Delta_U$ denotes the image of  the diagonal embedding $U \migiincl U \times U$.
 Also, 
for each $i=1,2$, we shall write  $\varpi_i$ for the projection $ (U \times U) \setminus \Delta_U \migi U$ to the $i$-th factor. 
The $\mcD_{\mbP^1 \times \mbA^1}$-module 
 $\msE := \eta_{\mbA^1 !*} (L\pi_1^* \int_\iota \msF \otimes^L L\mu^* \msK_0^\lambda)$ satisfies  
\begin{align} \label{Erfg}
\msE |_{\varpi_1^{-1}(U)} 
&\cong 
\left(\int_{\eta_{\mbA^1}} L\pi_1^* \int_\iota \msF \otimes^L L\mu^* \msK_0^\lambda \right) \Big|_{\varpi_1^{-1}(U)} \\
& \cong
  \int_{\eta_{\mbA^1}} \left(\int_{\iota_U}L\pi_{1, U}^* \msF \otimes^L L\mu^* \msK_0^\lambda \right) \Big|_{\varpi_1^{-1}(U)} \notag \\
  &\cong 
   \left(\int_{\kappa}\left(L\pi_{1, U}^* \msF \otimes^L L\mu^*\msK_0^\lambda \Big|_{U \times \mbA^1}\right) \Big|_{U \times U \setminus \Delta_U}\right)\Big|_{\varpi_1^{-1}(U)}, \notag
\end{align}
where 
\begin{itemize}
\item
$\iota_{U}$ denotes the open immersion $U \times \mbA^1 \migiincl \mbA^1 \times \mbA^1$ and  $\pi_{1, U}$ denotes the projection $U \times \mbA^1 \migi U$ onto the first factor, 
\item
the first ``$\cong$" follows from ~\cite[Corollary 2.8.5, (2)]{Kat2} (together with ~\cite[Theorem 7.1.1]{HTT}),
\item
the second ``$\cong$" follows from ~\cite[Theorem 1.7.3]{HTT}, and
\item 
the third ``$\cong$" follows from ~\cite[Corollary 1.7.5]{HTT}.
\end{itemize}
If we write $Z := (\mbP^1 \times \mbA^1) \setminus \mr{Im}(\kappa)$ (equipped with a structure of reduced scheme) and write $\zeta : Z \migiincl \mbP^1 \times \mbA^1$ for the natural closed immersion, then 
\begin{align}
\left(\int_\zeta \zeta^\dagger \msE \right) \Big|_{\varpi^{-1}_1 (U)}
& \cong \int_{\zeta} \zeta^\dagger (\msE |_{\varpi^{-1}_1 (U)}) \\
& \hspace{-1.5mm} \stackrel{(\ref{Erfg})}{\cong} 
\left(\int_\zeta \zeta^\dagger \left( \int_{\kappa}\left(L\pi_{1, U}^* \msF \otimes^L L\mu^* \msK_0^\lambda \Big|_{U \times \mbA^1} \right) \Big|_{U \times U \setminus \Delta_U}\right)
\right) \Big|_{\varpi_1^{-1}(U)} \notag \\
& = 0, \notag
\end{align}
where the last equality  follows from $\zeta^\dagger \int_\kappa (-) = 0$ asserted in   ~\cite[Proposition 1.7.1, (ii)]{HTT}.
In particular,   the equality   $\left(\int_{\pi_2} \int_\zeta \zeta^\dagger \msE\right) \Big|_{U} = 0$ holds.
Since the natural short exact sequence $0 \migi 
\int_\zeta \zeta^\dagger \msE \migi \msE \migi \int_\kappa \kappa^\dagger \msE \migi 0$
induces a distinguished triangle
\begin{align}
\left(\int_{\pi_2} \int_\zeta \zeta^\dagger \msE\right) \Big|_U
  \migi 
  \left(\int_{\pi_2}\msE\right) \Big|_U
   \migi
   \left( \int_{\pi_2}  \int_\kappa \kappa^\dagger \msE\right) \Big|_U \xrightarrow{+1},
\end{align}
$ \left(\int_{\pi_2}\msE\right) \Big|_U$ is quasi-isomorphic to
$\left( \int_{\pi_2}  \int_\kappa \kappa^\dagger \msE\right) \Big|_U$.
On the other hand,  we have  
 \begin{align} \label{HHH7}
\left(\int_{\pi_2} \int_\kappa \kappa^\dagger
\msE \right) \Big|_U
&\cong 
\int_{\varpi_2} \left(L\pi_1^* \int_\iota \msF \otimes^L L\mu^* \msK_0^\lambda\right) \Big|_{(U \times U) \setminus \Delta_U} \\
&\cong \int_{\varpi_2}  L\varpi_1^* \msF \otimes^L 
\left(L\mu^* \msK_0^\lambda |_{(U \times U) \setminus \Delta_U}\right). \notag 
  \end{align}
It follows that
the following sequence of  equalities holds:
\begin{align} \label{Ed88}
\rho (\mr{mc}_\lambda (\int_\iota \msF)) &= \rho (\int_{\pi_2} \msE) \\
&  = \rho (\left(\int_{\pi_2} \msE\right) \Big|_U) \notag \\
& =  \rho (\left( \int_{\pi_2}  \int_\kappa \kappa^\dagger \msE\right) \Big|_U) \notag \\
& = \rho (\int_{\varpi_2}  L\varpi_1^* \msF \otimes^L 
\left(L\mu^* \msK_0^\lambda |_{(U \times U) \setminus \Delta_U}\right)). \notag
\end{align}
  
Next,  by  Remark \ref{Rem4}, (iii),  and Proposition \ref{Prop558}, (ii), we  obtain  the inequalities 
 \begin{align}
 \rho_{\mbA^2_{\mcO_K}} (L\varpi_1^*\msF) & \leq \rho (\msF), \hspace{5mm} 
  \rho_{\mbA^2_{\mcO_K}} (L\mu^* \msK_0^\lambda |_{(U \times U) \setminus \Delta_U} )  \leq \rho_{} (\msK_0^\lambda) \leq H (\lambda), 
 \end{align}
 where $\mbA^2_{\mcO_K} := \mbA^1_{\mcO_K} \times_{\mcO_K} \mbA^1_{\mcO_K}$.
Hence, Remark \ref{Rem4}, (i), implies 
 \begin{align} \label{E46}
 \rho_{\mbA^2_{\mcO_K}} (L\varpi_1^* \msF \otimes^L 
L\mu^* \msK_0^\lambda |_{(U \times U) \setminus \Delta_U}) \leq  \rho (\msF) + H (\lambda).
 \end{align}

Also, since $L \varpi_1^* \msF \otimes^L 
 L\mu^* \msK_0^\lambda |_{(U \times U) \setminus \Delta_U}$ is of type $G$ (cf. Proposition \ref{Prop443}), it  has at most regular singularities and the exponents at any point at infinity on (the smooth proper model of) $U$ is contained in $\mbQ$
(cf. Remark \ref{Rem7}).
 It follows that, after possibly replacing $U$ with its open subscheme, 
 we can apply 
 ~\cite[Theorem 3.1.2]{AB} to the case where  ``$\mcM_{K_v}$" is taken to be the $v$-adic completion of the flat bundle $L \varpi_1^* \msF \otimes^L  
L \mu^* \msK_0^\lambda |_{(U \times U) \setminus \Delta_U} \left(=   \varpi_1^* \msF \otimes  
 \mu^* \msK_0^\lambda |_{(U \times U) \setminus \Delta_U}\right)$ for every $v \in \Sigma_{\mcO_K}$;
hence,  the inequality 
 \begin{align} \label{E45}
 \rho (\int_{\varpi_2}  L\varpi_1^* \msF \otimes^L 
L\mu^* \msK_0^\lambda |_{(U \times U) \setminus \Delta_U}) 
\leq 
(n^2 +1) \cdot  \rho_{\mbA^2_{\mcO_K}} (L \varpi_1^* \msF \otimes^L 
L \mu^* \msK_0^\lambda |_{(U \times U) \setminus \Delta_U})
 \end{align}
 holds.
 As a consequence,  
  the second inequality in (\ref{Er33}) can be verified by  combining (\ref{Ed88}), (\ref{E46}), and (\ref{E45}).
\end{proof}
\SSP

\begin{rema} \label{REM99}
If we know specific knowledge about a given chain complex $\msF^\bullet \in D_{h}^b (\mcD_{\mbA^1})_G$,
an upper bound for the generic inverse radius of its middle convolution may be obtained explicitly.
In fact, by induction on the cohomological length of $\msF^\bullet$,
the result of Proposition \ref{Pr34} enables us to reduce the situation to
 the case where $\msF^\bullet$ is a holonomic $\mcD_{\mbA^1}$-module $\mcF$.
Moreover, by induction on the length of the composition series of $\msF$,
it can be  assumed to be 
irreducible.
Since $\msF$ is the minimal extension  of a flat bundle on a locally closed smooth subvariety of $\mbA^1$,
$\msF$ is isomorphic to either $\int_q \msO$ (for some $q \in \mbA^1 (\overline{\mbQ})$) or
$\iota_{!*}\msG$ for some open immersion $\iota : U \migiincl \mbA^1$ and an irreducible  flat bundle $\msG$ on $U$.
 Thus, 
by applying  Proposition \ref{Prop558}, (ii),  and  Theorem  \ref{TT5},
we can estimate the value  $\rho (\mr{mc}_\lambda (\msF))$, as desired.
\end{rema}

\vspace{10mm}
\section{Equivalence among various  arithmetic properties on rigid flat bundles}  \label{S1133}
\LSP

The purpose of this final section 
is to prove  Theorem \ref{ThM}, asserting  a comparison among the classes of arithmetic flat bundles discussed so far.
We do so by restricting flat bundles  to rigid  ones  and applying Katz's middle convolution algorithm in order to reduce the problem to the rank one case.
We refer the reader to ~\cite{Ari} for a reasonable reference on  related topics.

Given a flat bundle $\msF$ on a nonempty open subscheme $U$ of $\mbA^1$ (where $\iota$ denotes the open immersion $U \migiincl \mbA^1$), 
we abuse notation by  writing  $\mr{mc}^0_\lambda (\msF)$ (where $\lambda \in \mbQ \setminus \mbZ$) for  the sheaf   $\mr{mc}_\lambda^0 (\int_\iota \msF) |_U$; this is none other than the classical definition of the middle convolution of  $\msF$.
In particular, it follows from a  well-known result  that $\mr{mc}^0_{- \lambda} (\mr{mc}^0_\lambda (\msF)) \cong \msF$.

\LSP
\subsection{Katz's middle convolution algorithm} \label{SS028}

Let $q$ be a closed point of the projective line $\mbP^1$ over $\overline{\mbQ}$.
Each  flat bundle $\msF$ on a nonempty open subscheme $U$ of $\mbP^1$ induces, via restriction,  
a flat bundle $\Psi_q (\msF)$ on  
the punctured formal neighborhood of $q$.
The isomorphism class  $[\Psi_q (\msF)]$ of $\Psi_q (\msF)$  is called
the {\bf formal type of $\msF$  at $q$}.
Then, we define 
the {\bf formal type} of $\msF$ to be  the collection
\begin{align}
\left\{ [\Psi_q (\msL)]\right\}_{q \in \mbP^1}.
\end{align}
Since $\Psi_q (\msF)$ is trivial when $q \in U$,
this collection is essentially determined by the subset $\left\{ [\Psi_q (\msF)]\right\}_{q \in \mbP^1\setminus U}$.

Given a flat bundle $\msF$ as above, we shall say that
 $\msF$ is {\bf rigid} (cf. ~\cite[Definition 2.2]{Ari}) if $\msF$ is determined by its formal type up to isomorphism, meaning that
any flat bundle  $\msF'$ on $U$ with  $\Psi_q (\msF) \cong \Psi_q (\msF')$ for every  $q \in \mbP^1$ is isomorphic to $\msF$.

We shall prove the following assertion based on the discussion in ~\cite[\S\,4]{Ari}.

\SSP
\bt \label{Th49}
Let $\msF$ be  an a.e.\,nilpotent flat bundle 
 on a nonempty  open subscheme $U \subseteq \mbA^1$.
Suppose that $\msF$ is irreducible and rigid, and that $\mr{rk}(\msF)>1$.
Then, 
there exists a pair 
\begin{align}
(\lambda, \msM)
\end{align}
 consisting of $\lambda \in \mbQ \setminus \mbZ$ and  a globally convergent flat bundle $\msM$ of rank one on $U$
such that
the flat bundle
$\mr{mc}_\lambda^0 (\msM \otimes \msF)$ 
has rank smaller than $\msF$.
\et
\begin{proof}
Let us take an arbitrary closed  point $q \in \mbP^1 \setminus U$.
Denote by  $\msV_q$ be the irreducible component $\msV'$ of 
the formal type $\Psi_q (\msF)$ of $\msF$ at $q$
 that minimizes the value 
\begin{align}
\frac{\delta (\mcH om (\msV', \Psi_q (\msF)))}{\mr{rk}(\msV')},
\end{align}
where $\delta (-)$ denotes the quantity defined in ~\cite[\S\,3.1]{Ari}.
Since $\Psi_q (\msF)$ is a.e.\,nilpotent,
it has at most regular singularities (cf. Remark \ref{Rem7}).
It follows that $\mr{rk}(\msV_q) =1$ and $\msV_q$ has regular singularities (cf. ~\cite[Example 4.1]{Ari}) and that   the residue  $\mr{res}(\msV_q)$ at $q$  lies in $\mbQ$.

Here, we shall suppose that the rational number  $\lambda := \sum_{q \in \mbP^1 \setminus U} \mr{res}(\msV_q)$ belongs to $\mbZ$.
Then, one can find a rank one flat bundle $\msN$ on $U$ with  $\Psi_q (\msN) \cong \msV_q$ for every $q \in \mbP^1 \setminus U$.
 By the Euler-Poincar\'{e} formula, either  $\mr{Hom}(\msF, \msN)$ or $\mr{Hom}(\msN, \msF)$ is nonzero (cf. ~\cite[Proposition 4.5]{Ari}).
Since $\msF$ is irreducible, this implies $\msF \cong \msN$.
In particular, 
$\msF$ has rank one, and this
 contradicts the assumptions.
It follows that $\lambda \in \mbQ \setminus \mbZ$.

Next, note that there exists a rank one flat bundle $\msM$ on $U$ satisfying 
\begin{align}
\Psi_q (\msM^\vee) \cong \begin{cases} \msV_q & \text{if $q \in \mbA^1 \setminus U$}, \\ 
\msV_\infty \otimes \hat{\msK}^{-\lambda}_\infty & \text{if $q = \infty$}, \end{cases}
\end{align} 
where $\hat{\msK}^{-\lambda}_\infty$ denotes the unique (up to isomorphism) regular singular  flat bundle of rank one 
on the punctured formal neighborhood of $\infty \in \mbP^1$
with $\mr{res}(\hat{\msK}^{-\lambda}_\infty) = -\lambda$.
Since $\mr{res}(\msV_q)$ lies in $\mbQ$ for each $q$, the flat bundle $\msM^\vee$ (hence also $\msM$) is globally convergent (cf. Proposition \ref{Prop443}).
Moreover, it follows from ~\cite[Propositions 3.6 and 4.6]{Ari} that
\begin{align}
\mr{rk}(\mr{mc}_\lambda (\mcH om (\msM^\vee, \msF))) \left(= \mr{rk} (\mr{mc}_\lambda (\msM \otimes \msF)) \right) < \mr{rk}(\msF).
\end{align}
This completes the proof of this proposition.
\end{proof}

\LSP
\subsection{Equivalence for  rigid Fuchsian systems} \label{SS0f28}

We shall conclude this paper by proving the following theorem.

\SSP
\bt[cf.  Theorem \ref{ThM}] \label{TT90}
Let  $U$ be  a nonempty  open subscheme of $\mbP^1$, and 
$\nabla$ a rigid connection on a vector bundle  over  $U$.
Then, the following three conditions are equivalent to each other:
\begin{itemize}
\item[(a)]
$\nabla$ is a $G$-connection;
\item[(b)]
$\nabla$ is a.e.\,nilpotent;
\item[(c)]
$\nabla$ is globally nilpotent.
\end{itemize}
\et
\begin{proof}
Since we already know the implications (a) $\Rightarrow$  (b) and (c) $\Rightarrow$  (b) (cf. (\ref{Eqww98})), it suffices to prove their inverse directions.
Let $\msF$ be an a.e.\,nilpotent rigid flat bundle on $U$.
We may assume, without loss of generality, that $\msF$ is irreducible; this is because the properties under consideration are all closed under taking flat subbundles, flat quotient bundles,  and  extensions of two flat bundles.
Also, after possibly shrinking $U$, we suppose that $U \subseteq \mbA^1$.

First, the case where $\msF$ has rank one follows from Proposition \ref{Prop562}.
Next, let us consider the case of $\mr{rk}(\msF) >1$.
By Theorem \ref{Th49},
there exists a sequence of flat bundles 
\begin{align} \label{Eq22j}
\msF = \msF_0 \mapsto \msF_1 \mapsto \cdots \mapsto \msF_n 
\end{align}
($n \geq 1$) on $U$ such that $\msF_{m+1}$ (for each $m=0, \cdots, n-1$) is obtained as  the result of the middle convolution
of $\msF_m$ with the parameter in $\mbQ \setminus \mbZ$ (i.e., $\msF_{m+1} := \mr{mc}^0_{\lambda_m} (\msF_m)$ for some $\lambda_m \in \mbQ \setminus \mbZ$) and that $\msF_n$ has rank one.
Since the operation of takings the middle convolution preserves the property of being a.e.\,nilpotent (cf. Theorem \ref{T491}),
$\msF_n$ turns out to  be a.e.\,nilpotent.
By the assertion for the rank one case (considered above), 
we see that $\msF_n$ is both globally nilpotent and of type $G$.
According to Theorem \ref{T491},  the inverse operation of each step in (\ref{Eq22j})  (which is described as middle convolution because $\mr{mc}^0_{-\lambda} (\mr{mc}_\lambda^0 (-)) \cong (-)$) preserves the property of being  globally nilpotent  (resp., of type $G$).
Thus,  $\msF$ is verified to be  globally nilpotent   (resp., of type $G$).
This completes the proof of this theorem.
\end{proof}

\LSP
\subsection*{Acknowledgements} \label{S04}

The author would like to  thank algebraic varieties over $\overline{\mbQ}$ for giving me a heartfelt encouragement and constructive comments about the globally inverse radius of a connection.
Our work was partially supported by Grant-in-Aid for Scientific Research (KAKENHI No. 21K13770).

\end{document}